\documentclass[amssymb,12pt ]{article}
\usepackage{geometry}
\usepackage{graphicx}
\usepackage{amsfonts}
\usepackage{amssymb}
\usepackage{amsmath}
\usepackage{longtable}

\usepackage{float}
\usepackage{amsbsy}
\usepackage{amstext}

\def\al {\aligned}
\def\eal {\endaligned}
\def\vtwo {\vskip 1pc}
\def\vv {\Vert}
\def\la {\langle}
\def\ra {\rangle}
\newtheorem{thm}{Theorem}[section]
\newtheorem{cor}[thm]{Corollary}
\newtheorem{lem}[thm]{Lemma}
\newtheorem{pro}[thm]{Proposition}

\newtheorem{defi}[thm]{Definition}
\newtheorem{con}[thm]{Convention}
\newtheorem{rem}[thm]{Remark}

\newenvironment{proof}{\noindent{\bf Proof. }}{\hfill\framebox[2mm]{}}

\begin{document}

\title{Trigonometry in extended hyperbolic space and extended de Sitter space}
\author{Yunhi Cho}
\date{}
\maketitle

\begin{abstract}
We study the hyperbolic cosine and sine laws in the extended
hyperbolic space which contains hyperbolic space as a subset and is
an analytic continuation of the hyperbolic space. And we also study
the spherical cosine and sine laws in the extended de Sitter space
which contains de Sitter Space $S^n_1$ as a subset and is also an
analytic continuation of de Sitter space. In fact, the extended
hyperbolic space and extended de Sitter space are the same space
only differ by $-1$ multiple in the metric. Hence these two extended
spaces clearly show and apparently explain that why many
corresponding formulas in hyperbolic and spherical space are very
similar each other. From these extended trigonometry laws, we can
give a coherent and geometrically simple explanation for the various
relations between the lengths and angles of hyperbolic polygons and
relations on de Sitter polygons which lie on $S^2_1$.
\end{abstract}

\renewcommand \figurename{Fig.}
\makeatletter
\renewcommand\fnum@figure[1]{\textbf{\figurename} }
\makeatother
\section{Introduction}

There are well known hyperbolic cosine and sine laws for triangles
in the hyperbolic space $\Bbb H^n$. If we consider Kleinian model
which contains the hyperbolic space as an open ball, we can think
about more general triangle which lies outside the hyperbolic space
or intersects the ideal boundary $\partial\Bbb H^n$. Then there is a
difficulty in geometric interpretation of such general type triangle
or other geometric objects. However the extended hyperbolic space
which is an analytic continuation of the hyperbolic space can give a
playground for such general geometric objects. Similarly extended de
Sitter space is obtained from de Sitter space $S^n_1$ and shows the
phenomena of the spherical geometry $\Bbb S^n$, just like the
extended hyperbolic space shows that of the hyperbolic geometry
$\Bbb H^n$.

In Section 2, we discuss what the extended model is and how it can
be constructed. The extended hyperbolic space which contains
hyperbolic space as a subset looks like the unit sphere $\Bbb S^n$
topologically. More detailed explanations about the extended space
are given in \cite{2}.

In Section 3, we explain how to define the notions of length and
angle on the extended space. In order to understand the extended
space more precisely, we should use $\epsilon$-approximation
technique. However here we only consider simple geometric objects
such as  length and angle, and we need not deeper theory of the
model. Here the length and angle must take complex values in
general. This kind of complex valued angle was introduced by Dzan
(\cite{3}, \cite{31}). He also constructed natural flat Lorentzian
geometry on $\Bbb R^{n,1}$ that looks like Euclidean geometry on
$\Bbb R^{n+1}$, then many formulas on $\Bbb R^{n,1}$ and $\Bbb
R^{n+1}$ exactly coincide each other. Schlenker \cite{7} also
defined complex valued distance and angle on Kleinian model using
 cross ratio. Our approach to distance and angle on the extended
space is more geometrically motivated and simple, and turned out to
be the same as Dzan and Schlenker's.

In Section 4, we prove the generalized hyperbolic (resp. spherical)
cosine and sine laws for the extended hyperbolic (resp. extended de
Sitter) space, those laws have exactly the same representation (see
Theorem \ref{3.10} and \ref{3.11}) of the original hyperbolic space
$\Bbb H^n$ (resp. spherical space $\Bbb S^n$). Note that most of the
proof and its difficulty for the generalized cosine and dual cosine
laws come from the sign ($\pm$) determining process.
 These generalized cosine and sine laws explain and easily deduce the well-known  relations (see Fenchel's book
\cite{4} or \cite{6} or \cite{8}) about the lengths and angles of
hyperbolic polygons in a simple unified way, for example, Lambert
quadrilateral, pentagon, rectangular hexagon, and so on. Furthermore
we can also obtain the similar relations between the lengths and
angles of de Sitter polygons on the pseudo-sphere (= Lorentz space
of constant curvature 1) $S^2_1$.

Lastly we remark some problems at the end of the paper which seem to
be important phenomena between the hyperbolic space $\Bbb H^n$ and
the spherical space $\Bbb S^n$.

\vskip 1pc {\noindent\bf Acknowledgement} The author would like to
thank to Hyuk Kim and Hyounggyu Choi. They gave some helpful
comments for this paper.

\section{Extended hyperbolic space and extended de Sitter space}

Our main concern is the unified trigonometry on the extended space,
so we should know what  the extended space is and why we need to
know the trigonometry on the model. For the answer of the first
question, the model is well explained in \cite{2} and reader can
easily understand the extended hyperbolic model itself and the
importance of the model. However we will introduce some necessary
parts of the theory in the following for convenience. And the second
question will be considered in Section 4.

To define and explain the extended model, let's start with the
hyperboloid model of hyperbolic space. Let $\Bbb R^{n,1}$ denote the
real vector space $\Bbb R^{n+1}$ equipped with the bilinear form of
signature $(n,1),$
$$
  \la x,y\ra =-x_0y_0+x_1y_1+ \cdots +x_ny_n,
$$
for all $x=(x_0,x_1,\cdots,x_n)$, $y=(y_0,y_1,\cdots,y_n)$.
Then the hyperbolic spaces $H^n_+$ and $H^n_-$, pseudo-sphere $S^n_1$ and light cone $L^n$ are defined by
$$
\al
 H^n_+&:= \{ x\in\Bbb R^{n,1} | \la x,x\ra =-1,\quad x_0>0\},\\
 H^n_-&:= \{ x\in\Bbb R^{n,1} | \la x,x\ra =-1,\quad x_0<0\},\\
 S^n_1&:= \{ x\in\Bbb R^{n,1} | \la x,x\ra =1\},\\
 L^n&:= \{ x\in\Bbb R^{n,1} | \la x,x\ra =0\}.
\eal
$$

We already know that $H^n_{\pm}$ has the induced Riemannian manifold
structure which has a constant sectional curvature $-1$, and that
$S^n_1$ becomes a Lorentzian  manifold (or semi-Riemannian of
signature ($-,+,\ldots ,+$))  which has a constant  sectional
curvature $1$, also called as de Sitter space (see \cite{5}). Now we
consider the Kleinian projective model. By the radial projection
$\pi_1$ with respect to the origin from $H^n_+$ onto $\{1\}\times
\Bbb R^n$, we obtain the induced Riemannian metric on the ball in
$\{1\}\times \Bbb R^n$  as follows (\cite{1},\cite{6}),
$$
ds^2_K=\left({\Sigma x_i dx_i\over 1-|x|^2}\right)^2 +{\Sigma dx_i^2\over 1-|x|^2}.
$$
If we extend this metric beyond  the unit ball using the same
formula, this metric  induces a semi-Riemannian structure outside
the unit ball in $\{1\}\times \Bbb R^n$. In fact, we compare this
metric with the one induced from the Lorentzian space  $S^n_1\cap
\{x=(x_0,x_1,\ldots ,x_n)|x_0>0\}$, by the radial projection into
$\{1\}\times \Bbb R^n$, then they differ only by sign $-1$. This
sign change of the metric implies the sign change of the sectional
curvature  from $+1$ to $-1$, which, of course, the curvature of the
metric $ds^2_K$. In this way,  we obtain an extended Kleinian model
with a singular metric $ds^2_K$ defined on $\{1\}\times \Bbb R^n$,
and this extended hyperbolic space $( \{1\}\times \Bbb R^n, ds^2_K)$
will be denoted by $K^n$.

In this paper, it is more convenient to consider the Euclidean unit
sphere in $\Bbb R^{n+1}$ with the induced metric coming from
$H^n_{\pm}$ and $-S^n_1$ ($S^n_1$ with $-1\times$ its metric) via
radial projection, and denote this model by $\Bbb S^n_H$. This {\it
hyperbolic sphere} model $\Bbb S^n_H$ on the Euclidean sphere
$\{x=(x_0,x_1,$ $\ldots ,x_n)|x_1^2+\cdots +x_n^2=1\}$ has three
parts: Two radial images of $H^n_{\pm}$, called the hyperbolic part,
as two open disks on upper and lower hemisphere and the radial image
of $S^n_1$, called the Lorentzian part, and these all three parts
have constant sectional curvature $-1$.

We can study the geometry of  $\Bbb S^n_H$ as an analytic
continuation of the hyperbolic space  $\Bbb H^n$. First  let's
define the volume of a set on the hyperbolic sphere. We denote
$dV_K$ and $dV_S$ for the volume forms on $K^n$ and $\Bbb S^n_H$
respectively. From the metric of the extended Kleinian model, we
have the following volume form  $dV_K$ (see \S6.1 of \cite{6}).

\begin{equation}\label{1}
\al
dV_K&=(\det(g_{ij}))^{\frac12}dx_1\wedge\cdots\wedge dx_n,\\
    &=\frac{dx_1\wedge\cdots\wedge dx_n}{(1-|x|^2)^{\frac{n+1}{2}}}.
\eal
\end{equation}

For any set $U$ on $\Bbb S^n_H\cap \{x\in \Bbb R^{n,1}|x_0>0\}$, we
can evaluate the volume of $U$ by
$$
\al
\text{vol} (U)&=\int_U dV_S\\
        &=\int_{\pi^{-1}(U)} dV_K\quad (\text{where }\pi \text{ is a radial projection: }K^n\to \Bbb S^n_H.)\\
        &=\int_{\pi^{-1}(U)}  \frac{dx_1\wedge\cdots\wedge dx_n}{(1-|x|^2)^{\frac{n+1}{2}}} \\
        &=\int_{G^{-1}(\pi^{-1}(U))} \frac{r^{n-1}}{(1-r^2)^{\frac{n+1}{2}}} dr\wedge d\theta,
\eal
$$
where $G:(r,\theta)\mapsto(x_1,\ldots ,x_n)$ is the polar
coordinates and  $d\theta$ is the volume form of the Euclidean
sphere $ \Bbb S^{n-1}$.

If $F(r)=\int_{G^{-1}(U)\cap S^{n-1}(r)}d\theta $ is an analytic
function of $r$, then the above integral becomes a 1-dimensional
integral as follows.
$$
 \int_{G^{-1}(U)} \frac{r^{n-1}}{(1-r^2)^{\frac{n+1}{2}}} dr d\theta=\int_a^b
\frac{r^{n-1}F(r)}{(1-r^2)^{\frac{n+1}{2}}} dr
$$
In general this integral does not make sense and the most natural
thing we can do instead is to define $\text{vol} (U)$ as the
following contour integral
\begin{con}\label{0}
\begin{equation}\label{f2}
\emph{vol}_H (U):=\int_{\gamma}
\frac{r^{n-1}F(r)}{(1-r^2)^{\frac{n+1}{2}}} dr
\end{equation}
\end{con}
where $\gamma$ is a contour from $a$ to $b$ for $a<1<b$ as depicted
below. Here we will fix its contour type as clockwise around $z=1$
once and for all throughout the paper.

\begin{figure}[h]
\begin{center}
\includegraphics[width=0.4\textwidth]{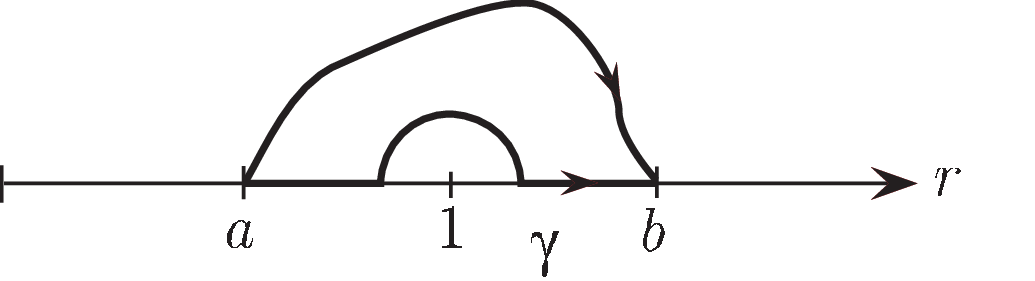}
\caption{\textbf{1}}
\end{center}
\end{figure}

Therefore we can compute the length of line segment on $\Bbb S^1_H$
by using the line integral (\ref{f2}). It easy to see that
\begin{equation}\label{3}
d_H(0,b):=\int^b_{0,\gamma} \frac{d r}{1-r^2}=\left\{\alignedat2
&\frac {1}{2}\log \frac{1+b}{1-b},\qquad &0\le b<1,\\
&\frac {1}{2}\log \frac{b+1}{b-1} +\frac{\pi}{2}i,\quad &1<b,
                           \endalignedat\right.
\end{equation}
where $d_H$ denotes 1-dimensional extended hyperbolic length of the
line segment $[0,b]$ in $K^1$ in the sense of (\ref{f2}).

If $F(r)$ is an analytic function around $r=1$, then it is easily
shown (see Proposition 2.1 of \cite{2}) that
$$
\lim_{\epsilon \to 0}\int_U dV_{K,\epsilon}:=\lim_{\epsilon \to
0}\int_U
\frac{d_{\epsilon}r^{n-1}}{(d_{\epsilon}^2-r^2)^{\frac{n+1}{2}}} dr
d\theta =\text{vol}_H (U), \quad\text{where }
d_{\epsilon}=1-\epsilon i~ (\epsilon>0).
$$
We called the above limit type approach $\epsilon$-approximation
technique which is more useful in the theoretical consideration. If
we choose $d_{\epsilon}=1+\epsilon i$ instead, then $\lim_{\epsilon
\to 0}\int_U d V_{K,\epsilon}$ will give a different value and a
slightly different geometry. That is exactly corresponding to a
contour integration with a {\it counterclockwise} around $z=1$,
i.e., going around at $z=1$ through lower half plane.

To determine the various geometric quantities which are to be
obtained as integrations on $\Bbb S^n_H$, the norms of vectors are
essential. From the sign change of the  metric on the pseudo-sphere
$S^n_1$, the norms of  tangent vectors on the Lorentzian part are
calculated by
$$\Vert x_p\Vert^2=-(-x_0^2+x_1^2+ \cdots +x_n^2),$$
and we should determine the sign of $\Vert x_p\Vert$ between plus and minus.
 On $\Bbb S^1_H$, we can determine  the sign by looking at the sign of  arc-length which can  be calculated
 by the contour integration formula (\ref{f2}) with $n=1$, i.e. $\int \frac{dr}{1-r^2}$,  and gives negative value outside $\Bbb H^1$.
 This gives us ($-1$) as the right choice of the sign of $\Vert x_p\Vert$ for the vectors in the radial direction on the Lorentzian
 part of
  $\Bbb S^n_H$. For the sign for the vectors normal to the radial direction, we use  the sign of 2-dimensional volume which is
  determined by one normal direction and one radial direction. It is not hard to check that the clockwise contour integral of the volume
  form gives sign $-i^{n-1}$ for Lorentzian part. Hence on the 2-dimensional spherical hyperbolic space $\Bbb S^2_H$, the volume for Lorentzian
  part has the sign $-i$ and thus the consistent choice of sign for the normal direction is $i$ (see Fig. 2).

\begin{figure}[h]
\begin{center}
\includegraphics[width=0.6\textwidth]{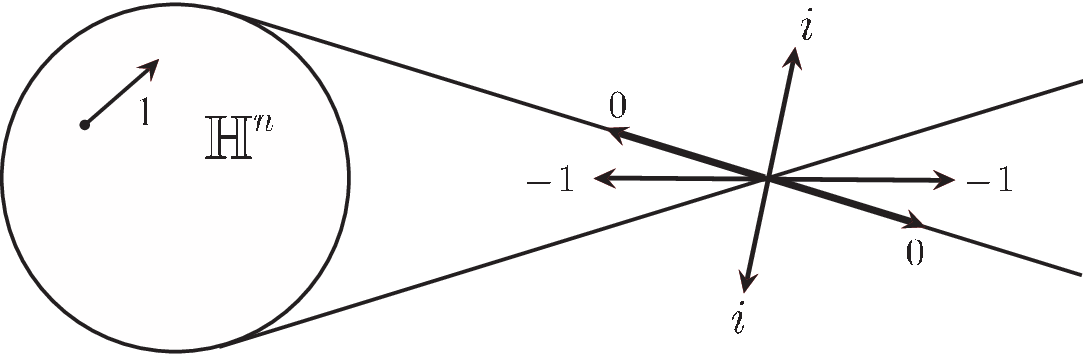}
\caption{\textbf{2}}
\end{center}
\end{figure}

\begin{con}\label{1.2} A tangent vector on the hyperbolic
part on $\Bbb S^n_H$  has a positive real norm, and a tangent vector
on the Lorentzian part on $\Bbb S^n_H$ has a negative real, zero, or
positive pure imaginary norm depending on whether it is timelike,
lightlike, or spacelike respectively.
\end{con}

Now let's think about another analytic continuation of the
pseudo-sphere $S^n_1$.   Basically the induced metric from
$S^n_1\cap \{x|x_0>0\}$ into  $\{1\}\times \Bbb R^n$ by the radial
projection differs by $-1$ from Kleinian metric $ds^2_K$, and we can
extend this metric to the inside of the unit ball. Note that we
always fix the sign of the norm of the tangent vector on $\Bbb
S^{n-1}=S^n_1\cap \{x|x_0=0\}$ as $+1$ as usual. Hence we will
assume the sign of the norm of the spacelike vector on $\Bbb S^n_S$
(see below) as $+1$. Here we denote the space and metric as $-K^n$
and $ds^2_{-K}$ respectively. Then by similar arguments we can
consider the unit sphere in $\Bbb R^{n+1}$ with induced metric
coming from the metric $ds^2_{-K}$ by the radial projection, and we
call this model as spherical sphere model and is denoted by $\Bbb
S^n_S$. Also we denote $dV_{-K}$ and $dV_{-S}$ as the volume forms
on $-K^n$ and $\Bbb S^n_S$ respectively, where $dV_{-K}=(-1)^{\frac
n2}dV_{K}$ and $dV_{-S}=(-1)^{\frac n2}dV_{S}$. Now we have to fix
the exact value of $(-1)^{\frac n2}$ between $i^n$ and $(-i)^n$. By
comparison of the norm of spacelike tangent vector at a Lorentzian
point on $\Bbb S^n_H$ and $\Bbb S^n_S$, $(-i)^n$ becomes a
reasonable choice between $i^n$ and $(-i)^n$.

\begin{con}\label{1.3} For any domain $U$ on $\Bbb S^n_S$, we
 evaluate the $n$-dimensional volume of $U$ by
$$
\emph{vol}_S(U)=\int_U dV_{-S}=(-i)^n \int_U dV_{S}=(-i)^n
\emph{vol}_H(U).
$$
\end{con}

From the similar chasing of volume form on $\Bbb S^1_S$, $\Bbb
S^2_S$  and considering of the equator of pseudo-sphere $S^n_1$,
i.e., the Euclidean  sphere $\Bbb S^{n-1}$, we naturally conclude
the following convention.

\begin{con}\label {1.4} A tangent vector on the hyperbolic
part on $\Bbb S^n_S$ has  a negative pure imaginary norm, and a
tangent vector on the Lorentzian part on $\Bbb S^n_S$ has a positive
pure imaginary, zero, or   positive real norm depending on whether
it is timelike, lightlike, or spacelike respectively (see Fig. 3).
\end{con}

\begin{figure}[h]
\begin{center}
\includegraphics[width=0.6\textwidth]{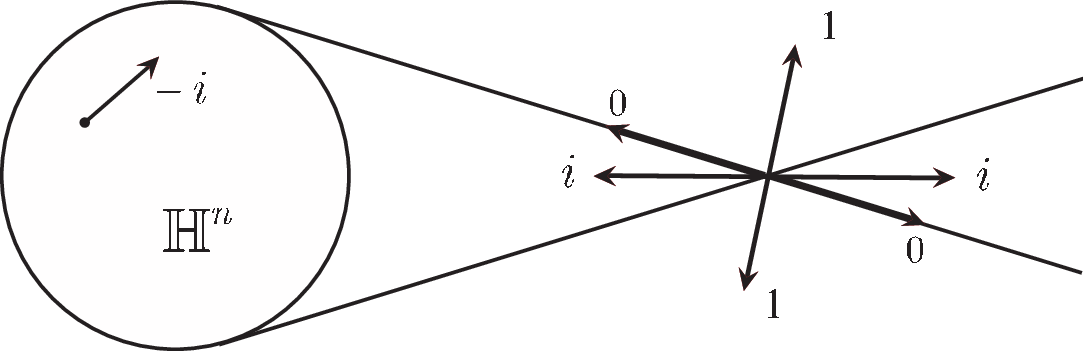}
\caption{\textbf{3}}
\end{center}
\end{figure}

We can see one of the similarities between $\Bbb S^n_H$, $\Bbb
S^n_S$ and $\Bbb S^n$ in the following theorem (see \cite{2} and
Convention \ref{1.3} for a proof). \vskip 1pc

\begin{thm}\label{1.5} $\emph{vol}_H(\Bbb
S^n_H$)$=i^n\cdot \emph{vol}(\Bbb S^n$) and $\emph{vol}_S(\Bbb
S^n_S)=\emph{vol}(\Bbb S^n$).
\end{thm}

If we change the contour type of the integral (\ref{f2}), we have
different relation between $\text{vol}_H(\Bbb S^n_H)$ and vol($\Bbb
S^n$).
 Also for the various different kinds of contour types, the conventions about $\Bbb S^n_S$ should be
 changed and the relations between $\text{vol}_S(\Bbb S^n)$ and vol($\Bbb S^n$) have similar modifications as  the hyperbolic ones too.
 If the contour is counterclockwise, then we have $\text{vol}_H(\Bbb S^n_H)=$$(-i)^n$ vol($\Bbb
 S^n$) (by slight change of the proof of Theorem 2.3 in \cite{2})
 and $\text{vol}_S(U)=i^n \text{vol}_H(U)$ (counterclockwise version of Convention \ref{1.3}).
 Hence we also get $\text{vol}_S(\Bbb S^n_S)=$vol($\Bbb S^n$).

For various kinds of contour types, we easily deduce the following
four formulas,
$$
\al
\text{vol}_H(\Bbb S^{2k-1}_H)&\equiv i^{2k-1} \text{vol }(\Bbb S^{2k-1}) \qquad (\text{mod }2 i^{2k-1}\text{ vol }(\Bbb S^{2k-1})),\\
\text{vol}_H(\Bbb S^{2k}_H)&=i^{2k} \text{vol }(\Bbb S^{2k}),\\
\text{vol}_S(\Bbb S^{2k-1}_S)&\equiv \text{vol }(\Bbb S^{2k-1}) \qquad (\text{mod }4 \text{ vol }(\Bbb S^{2k-1})),\\
\text{vol}_S(\Bbb S^{2k}_S)&=\text{vol }(\Bbb S^{2k}),\\
\eal
$$
Above formulas say that the total volume of even dimensional model
has unique value for any contour but odd dimensional model has
infinitely many values for various types of contours.

Note that we should know that two kinds of contours (clockwise and
counterclockwise contour) could be supported and be comprehended by
the $\epsilon$-approximation technique using $d_{\epsilon}=1\pm
\epsilon i$. But it is unclear that we can use an appropriate
$\epsilon$-approximation technique for other types of contours. So
we should make a proper mathematical theory to other contours.

Naturally the Lorentzian isometry group $O(n,1)$ can be considered
as the  isometry group of the hyperbolic sphere and spherical
sphere. More precisely, we know the following proposition (see
\cite{2} for a proof).

\begin{pro}\label{2.1} Let $U$ be a domain with piecewise analytic boundary transversal to $\partial\Bbb H^n$ in the
extended hyperbolic space. Then $\emph{vol}_H(U)$ has a well-defined
finite value and $\emph{vol}_H (g(U))=\emph{vol}_H (U)$ for each
$g\in PO(n,1)$.
\end{pro}

In fact, we already know that for a given $g$ in Isom($\Bbb H^n$),
which is  index two subgroup of $O(n,1)$, and for a given domain $U$
contained in $\Bbb H^n$, we get the equality $\text{vol}
(g(U))=\text{vol} (U)$. Surprisingly the volume of nice domains
intersecting with $\partial \Bbb H^n$($=\pi(L^n)$) can be
calculated. Though each part of the set divided by $\partial \Bbb
H^n$ has infinite volume, the total volume of two parts become
finite. This model has three infinite volume parts, $\pi(H^n_+)$,
$\pi(H^n_-)$, and $\pi(S^n_1)$, but by summing these parts we can
get a finite volume and hence a finite geometry without any
contradiction by using a finitely additive measure theory (see
\cite{2}).

\section{Length and angle on the extended
hyperbolic space and extended de Sitter space}

It is obvious from the definition of its metric that the geodesics
on $\Bbb S^n_H$  (resp. $\Bbb S^n_S$) are great circles on $\Bbb
S^n_H$ (resp. $\Bbb S^n_S$) and more generally the totally geodesic
subspaces are the intersections of the linear subspaces of $\Bbb
R^{n,1}$ with $\Bbb S^n_H$ (resp. $\Bbb S^n_S$) just like on $\Bbb
S^n$ (refer to \cite{5}).

 We denote  the distance between two points $A$ and $B$ in the
extended hyperbolic space $\Bbb S^n_H$ as $d_H(A,B)$. Let's first
discuss the distance between two points on $\Bbb S^1_H$. In this
case, the formula (\ref{3}) helps the calculation of the distance of
two points in  $\Bbb S^1_H$. For instance, if $A$ and $B$ are
symmetric with respect to the light cone $x_0=x_1$ in $\Bbb R^{1,1}$
as in Fig. 4 (i.e., $A$ and $B$ as vectors of $\Bbb R^{1,1}$ are
perpendicular), then their affine coordinates are $a(<1)$ and $\frac
1a$, and the distance will be $\frac{\pi}2i$ by formula (\ref{3}).
The distance between isometric images $A'$ and $B'$ of $A$ and $B$
will be again $\frac{\pi}2i$ being symmetric, and hence
$d_H(B,B')=-d_H(A,A')$ in Fig. 4.

\begin{figure}[h]
\begin{center}
\includegraphics[width=0.35\textwidth]{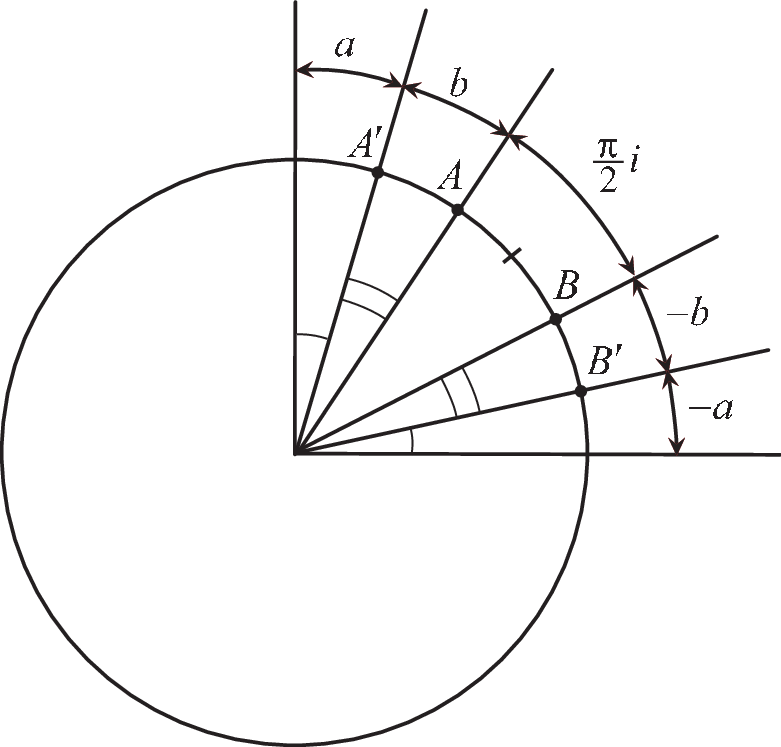}
\caption{\textbf{4}}
\end{center}
\end{figure}

To discuss the distance between two points in $\Bbb S^n_H$ in
general, it suffices to consider on $\Bbb S^2_H$.

For actual computations, it would be convenient to divide into the
following 3 cases. For the case when the geodesic connecting two
points meet $\partial\Bbb H^2$ transversely, we may assume that
these two points lie on $\Bbb S^1_H=\Bbb S^2_H\cap \{x|x_2=0\}$ by
an isometry and can handle as discussed above.

For the case when the geodesic line connecting these two points does
not intersect $\partial\Bbb H^2$, we can send this line to the
equator $(=\Bbb S^2_H\cap \{x|x_0=0\})$ of $\Bbb S^2_H$ by an
isometry, and hence the distance becomes $i$ times the distance on
the standard Euclidean unit circle.

The remaining case is when the line is tangent to $\partial\Bbb
H^2$. We can obtain the distance on the tangent line on $K^2$
through a theoretical way, but it needs a subtle
$\epsilon$-approximation technique (see \cite{2}). In this paper, we
consider the tangent case as a definition for  convenience.

\begin{defi}\label{3.1}For a point $x$ lie on $\partial \Bbb H^2$ and a dual geodesic $x^{\bot}$, the lengths of the
line segments in $x^{\bot}$ are defined by
$$
\al
d_H(w,y)&=0, \quad\text{if } w,y \text{ are in the same side with respect to }x,\\
d_H(y,z)&=\pi i, \quad\text{if } y,z \text{ are in the opposite sides with respect to }x,\\
d_H(x,y)&=d_H(x,z)=\frac{\pi}{2}i. \eal
$$
See Fig. 5.
\end{defi}

\begin{figure}[h]
\begin{center}
\includegraphics[width=0.5\textwidth]{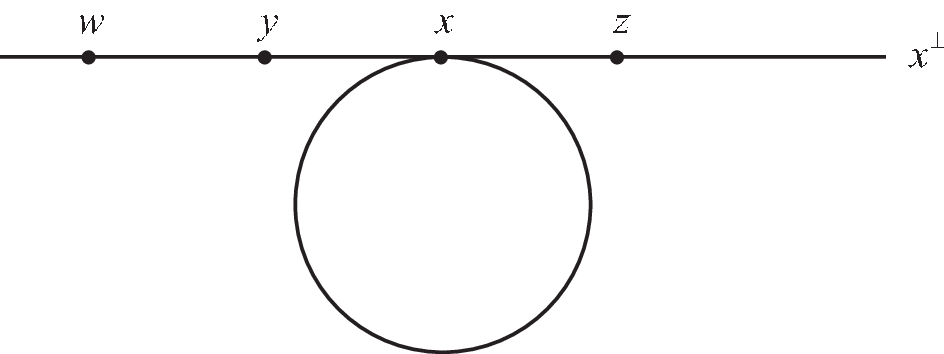}
\caption{\textbf{5}}
\end{center}
\end{figure}

In all of these discussions, we in fact have to show that
1-dimensional distance is invariant under isometry. That is shown in
Theorem 4.3  in \cite{2}.

We conclude the following theorem from the above discussion.
\begin{thm}\label{3.2} The total length of any great circle in $\Bbb
S^n_H$ (resp. $\Bbb S^n_S$) is $2 \pi i$ (resp. $2 \pi$).
\end{thm}

The extended hyperbolic space with Kleinian model has a projective
geometric structure, so a geodesic in the model is a straight line
and a dual of a point $x$, i.e., $x^{\bot}$  is easily obtained as
usual (see Fig. 6). Then the length of a geodesic line segment
joining $x$ (respectively $y$) and an arbitrary point  in $x^{\bot}$
(respectively $y^{\bot}$) is $\frac{\pi}2 i$. (Note if the model is
considered as a extended de Sitter space, then we should change
$\frac{\pi}2 i$ to $\frac{\pi}2$.) This follows since there is an
isometry which takes $x$ and $x^{\bot}$ to a point on the equator
and to a longitude respectively, and takes $y$ and $y^{\bot}$ to a
north pole and to the equator respectively.

\begin{figure}[h]
\begin{center}
\includegraphics[width=0.3\textwidth]{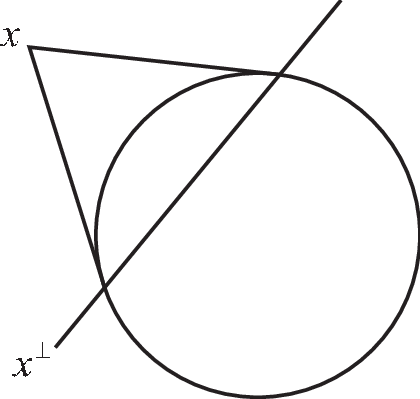}\quad\quad\quad
\includegraphics[width=0.35\textwidth]{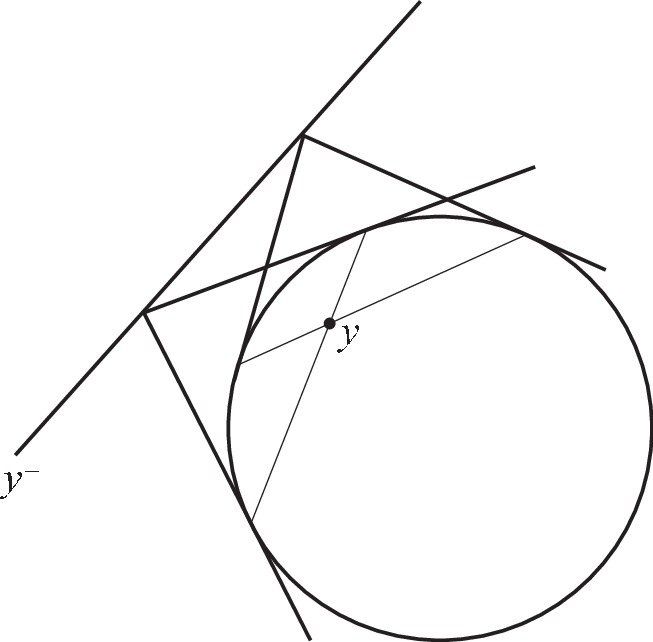}
\caption{\textbf{6}}
\end{center}
\end{figure}

Now we define angles on this extended model $\Bbb S^n_H$. From two
tangent vectors $v_p,w_p$  at a point $p$ on a Riemannian part, we
can define an angle $\theta$ by the equation

\begin{equation}\label{4} \la
v_p,w_p\ra= \vv v_p\vv \vv w_p\vv \cos\theta,\quad
0\le\theta\le\pi.\end{equation}

But for the Lorentzian part, we have some difficulties with this
formula since the function $\cos^{-1}$ is multi-valued and $\theta$
can take several complex values. The definitions of angle have been
given through the combinatorial way in \cite{3} and through the
cross ratio in \cite{7}.
  The following definition shows an
easy way of defining angle on $\Bbb S^n_H$ and $\Bbb S^n_S$. Note
that $v_p$ denotes the tangent vector at a point $p\in H^n_{\pm}$ or
$S^n_1\subset \Bbb R^{n,1}$ and $v\in \Bbb R^{n,1}$ independently
denotes the parallel translation of $v_p$ to origin.

\begin{defi}\label{2.2} For given two vectors $v,w \in \Bbb R^{n,1}$, the angle between $v$
and $w$, $\theta = \angle(v,w)$, is defined as $-i\cdot d_H(v,w)$
$(=d_S(v,w))$, where $d_H(v,w)$ (resp. $d_S(v,w)$) is the length of
a geodesic segment joining two points of $\Bbb S_H^n$ (resp. $\Bbb
S^n_S$) radially projected from $v,w$ to $\Bbb S_H^n$ (resp. $\Bbb
S^n_S$).

\begin{figure}[h]
\begin{center}
\includegraphics[width=0.5\textwidth]{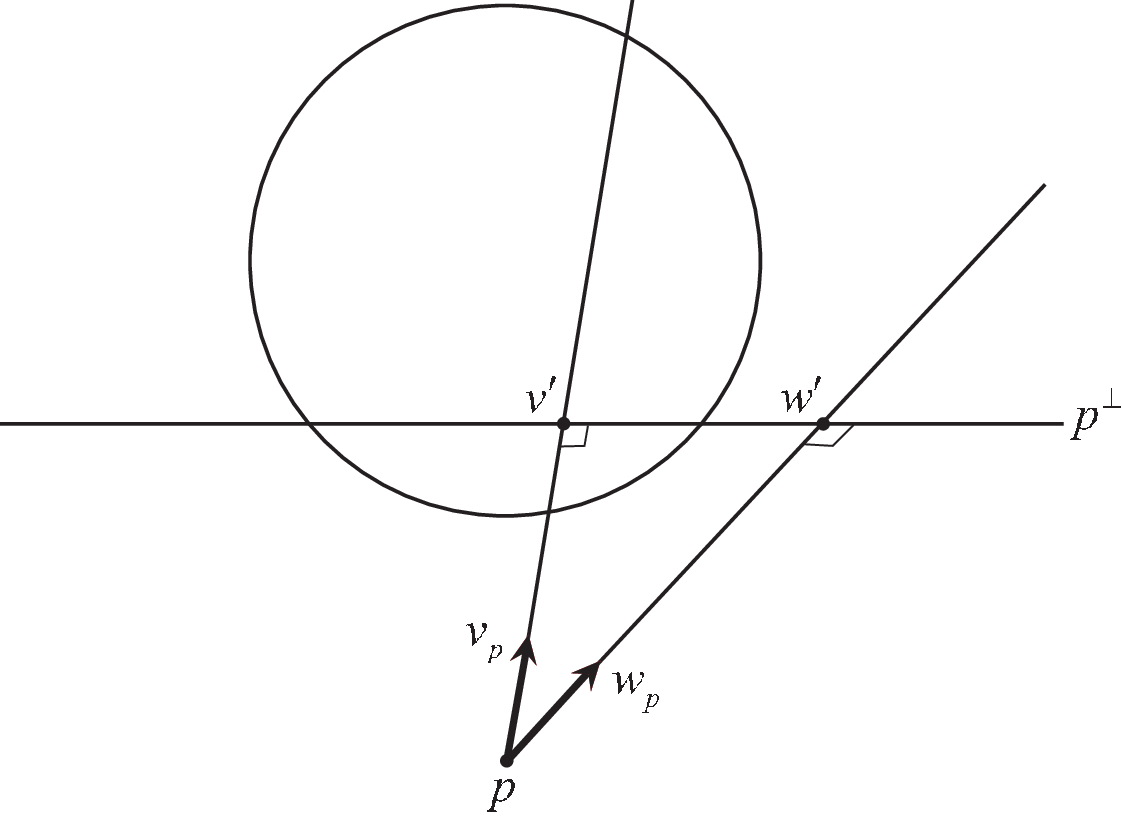}
\caption{\textbf{7}}
\end{center}
\end{figure}

For $p\in \Bbb S^n_H$ and $p\notin\partial\Bbb H^n$, the angle
$\angle(v_p,w_p)$ between two tangent vectors $v_p,w_p\in T_p\Bbb
S^n_H$ is defined as $-i\cdot d_H(v',w')$, where $v'$ is a point
which is obtained by the intersection of the dual plane $p^{\bot}$
and the geodesic line starting at $p$ with direction $v_p$ (see Fig.
7).

\begin{figure}[h]
\begin{center}
\includegraphics[width=0.35\textwidth]{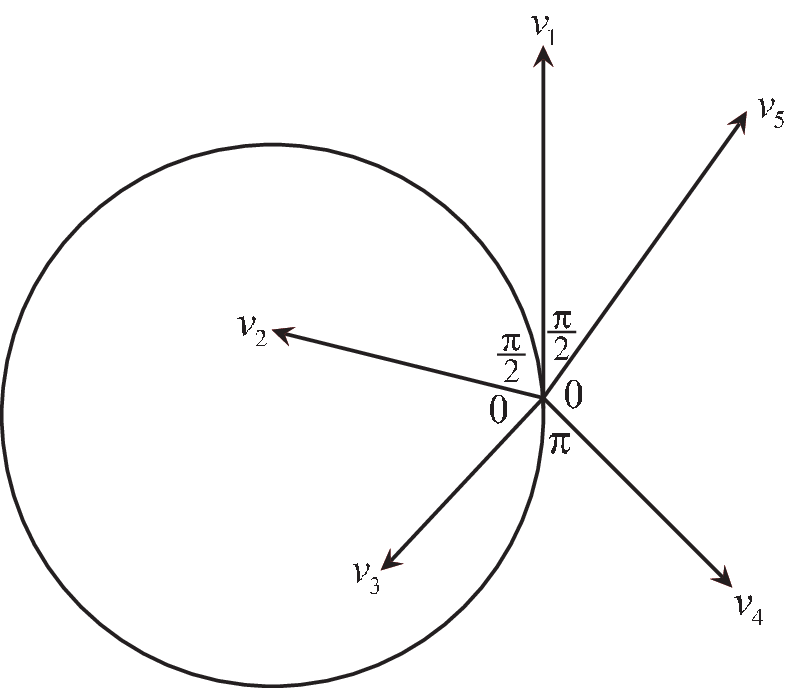}
\caption{\textbf{8}}
\end{center}
\end{figure}

If $p\in\Bbb S^2_H$ lies on $\partial\Bbb H^2$, we define the angle
$\angle(v_i,v_{i+1})$ for the tangent vectors $v_1,\ldots,v_5$
configured as in Fig. 8 as follows:

\noindent $$\angle(v_1,v_2)=\angle(v_5,v_1)=\frac {\pi}{2},~~
\angle(v_2,v_3)=\angle(v_4,v_5)=0, \text{ and
}\angle(v_3,v_4)=\pi.$$

If $p \in \partial \Bbb H^n$ with $n \geq 3 $, we have another type
of tangent plane at $p$ which touches $\partial \Bbb H^n $ at the
only point $p$. In this plane, we define the angle at $p$ as the
usual Euclidean angle.
\end{defi}

It is clear that an angle of one rotation around a point is $2\pi$
by Definition \ref{2.2} and Theorem \ref{3.2}. Notice that the
factor $-i$ is multiplied to normalize the total length $2\pi i$ of
the great circle in $\Bbb S_H^n$ as $2 \pi$ (see Theorem \ref{3.2}).

The isometry invariance of an angle at a point $p\in \Bbb
S^n_H\backslash \partial\Bbb H^n$ is obtained from the invariance of
distance.

\begin{rem}\label{mm}In fact, the second part of the definition is obtained from the first part of the definition,
but we made it as a definition for  convenience. Even though a
justification of the third and fourth part of the definition comes
from the $\epsilon$-approximation technique, we only refer the
reader to \cite{2}.
\end{rem}

\begin{defi}\label{2.3}
 For a lune $l(x_p,y_p)$, $\angle
(x_p,y_p)$ denotes the angle of $l(x_p,y_p)$ at the vertex $p$. Here
the lune $l(x_p,y_p)$ is the inner region generated by two half
great circles starting at $p$ with direction $x_p,y_p$ respectively
and ending at $-p$.

A lens $L(x^{\bot},y^{\bot})$ is the intersection of two hemispheres
$H_x$ and $H_y$,  where the hemisphere $H_x$ is posed opposite to
$x$ and $\partial H_x$ is perpendicular to $x$, and $\angle
(x^{\bot},y^{\bot})$ denotes the dihedral angle of the lens
$L(x^{\bot},y^{\bot})$. A lens $L(x^{\bot},y^{\bot})$ is called
ideal if $H_x\cap H_y$ meets the $\partial\Bbb H^n$ at two points
only.
\end{defi}

\begin{defi}  For a given lens $L(x^{\bot},y^{\bot})$, the dihedral angle of the lens is defined as the
angle $\angle (u_p,v_p)$. Here $u_p$ (resp. $v_p$) is a tangent
vector on $\partial H_x$ (resp. $\partial H_y$) with a base point
$p\in H_x\cap H_y$, and $u_p, v_p$ are perpendicular to $H_x\cap
H_y$. Note that for a non-ideal lens case the vertex $p$ can take
any point in $H_x\cap H_y$, for an ideal lens case the vertex $p$
only can take one of two ideal points in $H_x\cap H_y$.
\end{defi}

\begin{rem} It is easy to show the well-definedness of the dihedral angle of a lens. In particular for a non-ideal
lens, tangent vectors $u_p,v_p$ are uniquely determined up to
positive constant magnitude. But for an ideal lens, tangent vectors
$u_p,v_p$ at an ideal point $p$ can have infinitely many directions.

Since the dihedral angle of a polyhedron can be defined as the
dihedral angle of the induced lens naturally. The notion of dihedral
angle becomes an important object in the polyhedron theory at the
extended hyperbolic space.
\end{rem}

A lune is a 2-dimensional object and a lens in $\Bbb S^n_H$ or $\Bbb
S^n_S$ is  an $n$-dimensional object.

If a lune $l(x_p,y_p)$ with an angle $\theta$ and a lens
$L(x^{\bot},y^{\bot})$  with a dihedral angle $\alpha$ are given in
$\Bbb S^2_H$ or $\Bbb S^2_S$, then by case by case examinations we
get one of the following three kinds of relations (see Fig. 9):
$$
\alpha =\pi-\theta \text{ or }-\pi+\theta \text{ or }\pi+\theta.
$$

\begin{figure}[h]
\begin{center}
\includegraphics[width=0.6\textwidth]{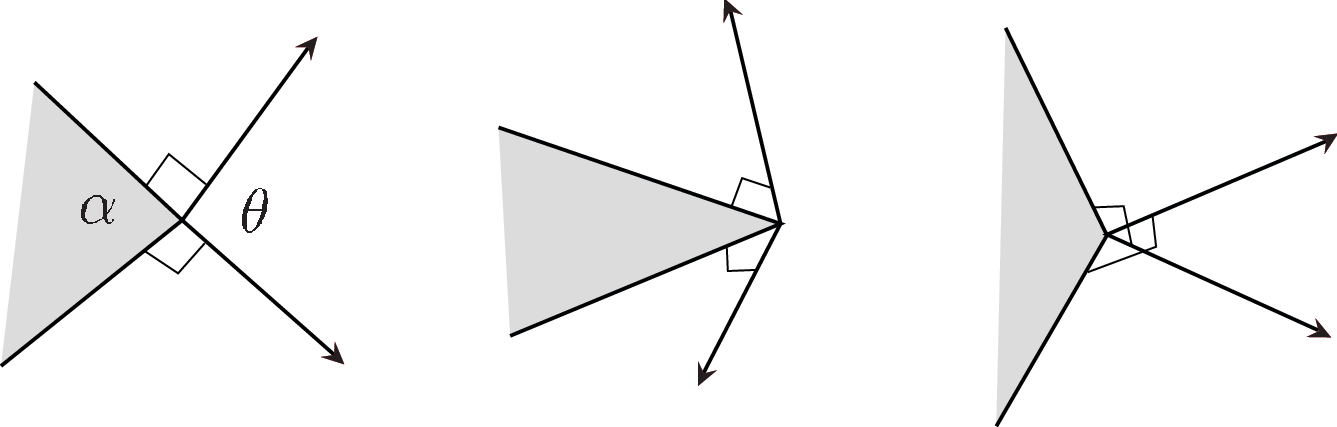}
\caption{\textbf{9}}
\end{center}
\end{figure}

 We already know that the Riemannian case has unique
relation $\alpha =\pi-\theta$.  For an $n$-dimensional lens, we also
conclude the same result as the 2 dimensional case.

\begin{lem}\label{2.4} For a lens $L(x^{\bot},y^{\bot})$,
there are equalities, $$
\angle(x^{\bot},y^{\bot})=\pm(\pi-\angle(x_{p},y_{p}))\text{ or
}\pi+ \angle(x_{p},y_{p}).$$
\end{lem}

\begin{proof} The dihedral angle of non-ideal lens $L(x^{\bot},y^{\bot})$ is the
same as  an angle of a lune which is the intersection of the lens
$L(x^{\bot},y^{\bot})$ and a 2-dimensional embedded geodesic
plane(in fact, a 2-sphere) generated by a point $p$ and two tangent
vectors $x_{p},y_{p}$, where the point $p$ is an arbitrary point in
$x^{\bot}\cap y^{\bot}$. So we can apply the 2-dimensional result to
the $n$-dimensional case.

For an ideal lens $L(x^{\bot},y^{\bot})$, we can check the relation
$\angle(x^{\bot},y^{\bot})=\pi-\angle(x_{p},y_{p})$.
\end{proof}

\begin{rem}\label{choi} For a given lens with dihedral angle
$\alpha=\angle(x^{\bot},y^{\bot})$, if we define an oriented angle
$\tilde\theta$ and can pass over the $\pi$,  then we can unify the
three relations to a single relation $\alpha =\pi-\tilde\theta$.
\end{rem}

The following lemma is given in Thurston's book \cite{8} or
\cite{6}.
 At first we need some notations: The Lorentzian norm of a vector $x$ in $\Bbb R^{n,1}$ is defined to be a complex number
$$ \Vert x\Vert={\langle x,x \rangle}^{1/2},$$
where $\vv x\vv$ is either positive, zero, or positive pure
imaginary. If $\Vert x \vv$ is positive imaginary,  we denote its
absolute value by $\vert\vv x\vv\vert$. In fact, $ \Vert x\Vert$ can
have minus or minus pure imaginary values, but those choices are not
suitable to our clockwise contour convention \ref{0}.

 We have to
be cautious about the difference between $\vv x\vv$ and $\vv x_p\vv$
for $x\in \Bbb R^{n,1}$. The vector $x$ is parallel translation of
the tangent vector $x_p\in T_p S^n_1$ or $T_p H^n_{\pm}$ to the
origin. If a point $p$ is contained in the Lorentzian part of $\Bbb
S^n_H$ (resp. $\Bbb S^n_S$), then we know $\vv x_p\vv =i \vv x\vv$
(resp. $\vv x_p\vv =\vv x\vv$) by Convention \ref{1.2} and
\ref{1.3}. Since we changed the sign of the induced metric on the
Lorentzian part. Also  if a point $p$ is contained in the hyperbolic
part of $\Bbb S^n_H$ (resp. $\Bbb S^n_S$), then we know $\vv x_p\vv
=\vv x\vv$ (resp. $\vv x_p\vv =-i \vv x\vv$) by Convention \ref{1.2}
and \ref{1.3}. Note the sign change of metric induces $\langle
x_p,y_p \rangle = -\langle x,y \rangle$. Hence we have the identity:
$$\frac{\la x_p,y_p\ra}{\vv x_p\vv \vv y_p\vv}=\frac{\la
x,y\ra}{\vv x\vv \vv y\vv}.$$

\begin{lem}\label{2.6} (interpretation of the inner product) If
$x$ and $y$ are vectors in $\Bbb R^{n,1}$, then either
\newline
(a) $x,y$ are timelike vectors and   $\la x,y\ra =\pm\vv x\vv \vv
y\vv\cosh d_H(\pm x,y)$; or
\newline
(b) $x$ is a timelike vector and $y$ is a spacelike vector, and

$ \la x,y\ra =\pm \vert\vv x\vv\vert\ \vv y\vv\sinh d_H(x,y^{\bot})$; or
\newline
(c) $x,y$ are spacelike vectors and the hyperplanes
$x^{\bot},y^{\bot}$ are secant, parallel

or ultra parallel depending on whether the intersection $x^{\bot}\cap\ \  y^{\bot}$ pass

through respectively the inside  of $\Bbb H^n, \partial \Bbb H^n$ or the outside of $\Bbb H^n$ only. In

the first case, $\la x,y\ra \ =\ -\vv x\vv \vv y\vv\ \cos \angle(x^{\bot},y^{\bot})$; in the second, $\la x,y\ra =$

$\pm\vv x\vv \vv y\vv$; and in the third, $\la x,y\ra =\pm\vv x\vv \vv y\vv\cosh d_H(x^{\bot},y^{\bot})$.
\end{lem}

The above Lemma has many cases for explaining the inner product.
However our new notion $d_H(x,y)$ in the extended hyperbolic space
enables us unify all these cases into a single form as in the
following theorem. This shows one good aspect of natural property of
the extended space.

\begin{thm}\label{2.10} For vectors $x$ and $y$ in
the Lorentzian space $\Bbb R^{n,1}$ and with condition
$d_H(x,y)\ne\infty$, we have

$$\la x,y\ra = \vv x\vv \vv y\vv\cosh
d_H(x,y)$$ Note that the case $d_H(x,y)=\infty$ induces that $\cosh
d_H(x,y)=\infty$, $\vv x\vv$ or $\vv y\vv=0$ (one of two vectors $x$
and $y$ becomes a lightlike vector), and $\la x,y\ra$ take a certain
value, hence we get $0\cdot\infty=\text{constant}$. In some sense,
the formula is always true for all cases.
\end{thm}

\begin{proof} It is sufficient to consider the following 3 cases by
isometric changes. The cases are respectively when the intersection
of the plane $span \{v,w\}$ and $\Bbb S^n_H$ is $\Bbb S^1_H$, the
equator $\Bbb S^2_H\cap \{x|x_0=0\}$, or the great circle tangent to
$\partial\Bbb H^2$.

By linear property of $\la\cdot,\cdot\ra$ and $\vv\cdot\vv$, and the
definition of $d_H$, we can assume that $\vv x\vv, \vv y\vv=1$ or
$i$ or $0$.

For the first case, let's suppose $x=(1, 0)$, then if $y$ is
timelike vector, then $y$ is represented by $(\pm\cosh a, \sinh a)$;
if $y$ is spacelike vector, then $y$ is represented by $(\pm\sinh a,
\cosh a)$. Also we should consider two spacelike vectors case $x=(0,
1)$ and $y=(\sinh a, \pm\cosh a)$.

Second case induces $x=(0, 1, 0)$ and $y=(0, \cos a, \sin a)$.

Third case induces $x=(1, 1, 0), y=(a, a, 1)$ or $x=(1, 1, 0),
y=(-1, -1, 0)$ or $x=(a, a, 1), y=(b, b, 1)$ or $x=(a, a, 1), y=(b,
b, -1)$.  Note all above $a, b$ are positive numbers.

All cases are checked  below.

\begin{itemize}
\item $x=(1, 0)$, $y=(\cosh a, \sinh a)$ implies $d_H(x, y)=d_H(0, \tanh a)=\log \sqrt{\frac{1+\tanh a}{1-\tanh
a}}$ hence $\vv x\vv \vv y\vv\cosh d_H(x,y)=i\cdot i\cdot \cosh
a=-\cosh a=\la x,y\ra$.
\item $x=(1, 0)$, $y=(-\cosh a, \sinh a)$ implies $d_H(x, y)=\pi i-d_H(0, \tanh a)$ hence $\vv x\vv \vv y\vv$ $\cosh d_H(x,y)$ $=i\cdot i\cdot (-\cosh
a)=\cosh a=\la x,y\ra$.
\item $x=(1, 0)$, $y=(\pm\sinh a, \cosh a)$ implies $d_H(x, y)=\frac{\pi}{2} i\pm d_H(0, \tanh a)$ hence $\vv x\vv \vv y\vv$ $\cosh d_H(x,y)$ $=i\cdot 1
\cdot (\pm i\sinh a)=\mp\sinh a=\la x,y\ra$.
\item $x=(0, 1)$, $y=(\sinh a, \cosh a)$ implies $d_H(x, y)=-d_H(0, \tanh a)$ hence $\vv x\vv \vv y\vv$ $\cosh$ $d_H(x,y)$ $=1\cdot 1
\cdot \cosh a=\cosh a=\la x,y\ra$.
\item $x=(0, 1)$, $y=(\sinh a, -\cosh a)$ implies $d_H(x, y)=\pi i+d_H(0, \tanh a)$ hence $\vv x\vv \vv y\vv$ $\cosh$ $d_H(x,y)$ $=1\cdot 1
\cdot (-\cosh a)=-\cosh a=\la x,y\ra$.
\item $x=(0, 1, 0)$, $y=(0, \cos a, \sin a)$ implies $d_H(x, y)=ai$ hence $\vv x\vv \vv y\vv$ $\cosh$ $d_H(x,y)$ $=1\cdot 1
\cdot \cos a=\cos a=\la x,y\ra$.
\item $x=(1, 1, 0)$, $y=(a, a, 1)$ implies $d_H(x, y)=\frac{\pi}{2}i$ hence $\vv x\vv \vv y\vv$ $\cosh$ $d_H(x,y)$ $=0\cdot 1
\cdot \cosh \frac{\pi}{2}i=0=\la x,y\ra$.
\item $x=(1, 1, 0)$, $y=(-1, -1, 0)$ implies $d_H(x, y)=\pi i$ hence $\vv x\vv \vv y\vv$ $\cosh$ $d_H(x,y)$ $=0\cdot 0
\cdot \cosh \pi i=0=\la x,y\ra$.
\item $x=(a, a, 1)$, $y=(b, b, 1)$ implies $d_H(x, y)=0$ hence $\vv x\vv \vv y\vv$ $\cosh$ $d_H(x,y)$ $=1\cdot 1
\cdot \cosh 0=1=\la x,y\ra$.
\item $x=(a, a, 1)$, $y=(b, b, -1)$  implies $d_H(x, y)=\pi i$ hence $\vv x\vv \vv y\vv$ $\cosh$ $d_H(x,y)$ $=1\cdot 1
\cdot \cosh \pi i=-1=\la x,y\ra$.
\end{itemize}

Now we have examined all the cases and complete the proof.
\end{proof}

Here we do not use the result of Lemma \ref{2.6} in the proof of
Theorem \ref{2.10}. However we can prove Theorem \ref{2.10} from
Lemma \ref{2.6}. Conversely, we can prove Lemma \ref{2.6} from
Theorem \ref{2.10}.

\begin{cor}\label{2.11} For vectors $x$ and $y$
in the Lorentzian  space $\Bbb R^{n,1}$ and a point $p$ in
$x^{\bot}\cap y^{\bot}$ and with condition $d_H(x,y)\ne\infty$, we
have
 \begin{equation}\label{5}\la x,y\ra=\vv x\vv \vv y\vv\cos
\angle(x,y),\end{equation}
$$
\aligned
 \la x_p,y_p\ra &= \vv x_p\vv \vv y_p\vv\cos \angle(x_p,y_p),\\
 \la x,y\ra &= -\vv x\vv \vv y\vv\cos \angle(x^{\bot},y^{\bot}),\\
 \la x,y\ra &= \vv x\vv \vv y\vv\cos d_S(x,y).
\endaligned
$$
\end{cor}

\begin{proof} See Lemma \ref{2.4} and Definition \ref{2.2}.
\end{proof}

Corollary \ref{2.11} shows that the hyperbolic sphere $\Bbb S^n_H$
and the spherical sphere $\Bbb S^n_S$ and the definitions about
length and angle on the spaces have natural and essential
properties.

We already showed $\la x,y\ra=\vv x\vv \vv y\vv\cos \angle(x,y)$
from Definition \ref{2.2}. If we add the following three properties
to the formula (\ref{5}), then we can show that the angle is
uniquely determined by these four properties and equivalent to
Definition \ref{2.2}. The additional three properties are

\noindent (i) the invariance under isometry,\\
\noindent (ii) finite additivity of angle: if $\theta$ consists of
two parts $\theta_1$ and $\theta_2$, then
$\theta=\theta_1+\theta_2$,\\
\noindent (iii) the angle of half rotation is $\pi$, i.e., a
straight line has angle $\pi$.

Other equivalent angle definitions are shown at Remark 4.13 in
\cite{2}.

\section{Cosine laws and sine law for general
triangles}

\leftline{\bf 4.1 Cosine laws} \vtwo
In the hyperbolic space, cosine
laws and sine law are basic laws as well as the spherical space.  So
we have to examine whether these laws are satisfied in the extended
hyperbolic space.

We need some definitions.

\begin{defi}\label{4.1} For a negative real number $a$,
$\sqrt{a}$ means $\sqrt{-a} i$.
\end{defi}\label{}
For example, $\sqrt{4}=2$ and $\sqrt{-4}=2i$.
\begin{defi} For a complex number $a\in(\Bbb R\cup\Bbb R i)-\{0\}$,
$\emph{sgn }(a)$ is defined as follows:
$$\emph{sgn } (a)=\left\{\alignedat2
1 \quad&\text{if $a$ is positive or positive pure imaginary},\\
-1 \quad&\text{if $a$ is negative or negative pure imaginary}.
                           \endalignedat\right.$$
\end{defi}
The sgn notation is slightly generalized, so the usual properties
are not satisfied any more. For example,
$\text{sgn}(ab)=\text{sgn}(a)~\text{sgn}(b)$ is not satisfied, if
both of $a$ and $b$ are pure imaginary numbers.
\begin{defi}\label{4.2} For non-zero real numbers
$a_1,a_2,\ldots,a_n$, the function $\emph{msgn}$(many elements sign)
is defined by
$$ \emph{msgn}(a_1,a_2,\ldots,a_n)=\frac{\sqrt{a_1}\sqrt{a_2}\cdots\sqrt{a_n}} {\sqrt{a_1a_2\cdots a_n}}.$$
\end{defi}

From the Definition \ref{4.2}, we easily obtain the following
proposition.
\begin{pro}\label{3.3} For  non-zero real numbers $a,a_1,\ldots,a_n,b_1,\ldots,b_m$,
we obtain
\newline
(a) $ \emph{msgn}(a)=1,$
\newline
(b) $ \emph{msgn}(a,a)=\emph{sgn}(a),$ and
$\emph{msgn}(a_1,a_2,\ldots,a_n)~ \emph{msgn}(a_1,a_2,\ldots,a_n)=1$
\newline
(c) $ \emph{msgn}(a_1,a_2,\ldots,a_n)\
\emph{msgn}(b_1,b_2,\ldots,b_m)$

$= \emph{msgn}(a_1,\ldots,a_n,b_1,\ldots,b_m)\ $$
\emph{msgn}(a_1\cdots a_n,b_1\cdots b_m),$
\newline
(d) $ \emph{msgn}(a_1,a_2,\ldots,a_n)=(-1)^{[\frac{\alpha}2]}$,
where $[\cdot]$ is the Gauss notation and $\alpha$ is the number

of negative elements among $a_i, i=1,\ldots,n$,
\newline
 (e) $\emph{msgn}(a_1,a_1,a_2,a_2,\ldots,a_n,a_n)=\emph{sgn}(a_1a_2\cdots
 a_n)$.
\end{pro}
\begin{proof} All of these follows easily from Definition \ref{4.2}.
\end{proof}

We also need next definitions to prove the cosine and sine laws.

\begin{defi}\label{3.4}\
\newline (1) For a given hemisphere $H\in \Bbb S^n_H$ (resp. vector $v\in
\Bbb S^n_H$), the algebraic dual

of $H$ (resp. $v$) is a point $v$ (resp. hemisphere $H$) given by $
\la H,v\ra\ge 0 $, i.e., ${}^{\forall} h\in H$

$ \la h,v\ra\ge 0$.
\newline
(2) For a given hemisphere $H\in \Bbb S^n_H$ (resp. vector $v\in \Bbb S^n_H$), the geometric

 dual of $H$ (resp. $v$) is a point $v$ (resp. hemisphere $H$) given by $\la \partial H,v\ra=0$

and $v\notin H$.
\end{defi}
For convenience sake, we denote an algebraic (resp. geometric) dual
of $X$ as $X^{a\bot}$ (resp. $X^{g\bot}$).

\begin{rem}
For a given ideal (i.e. tangent to $\partial\Bbb H^n$) hemisphere
$H$, the algebraic dual vector $v=H^{a\bot}$ is well defined. But
the geometric dual vector $v=H^{g\bot}$ is not well defined and
there are two direction choices (i.e. two points in  $\Bbb S^n_H$)
for the vector $v$.
\end{rem}

\begin{defi}\label{3.5} For given three linearly independent
non-lightlike vectors $v_1,v_2,v_3$  (resp. non-ideal hemispheres
$H_1,H_2,H_3$) in $\Bbb S^2_H$, the three vectors (resp.
hemispheres) induce a unique triangle $\triangle(v_1,v_2,v_3)$
(resp. $\triangle(H_1,H_2,H_3)$) with sides composed of
 \textsl{"smaller"} geodesics (resp.
$H_1\cap H_2\cap H_3$). Then the algebraic (resp. geometric) dual of
triangle $\triangle=\triangle(H_1,H_2,H_3)$ is a triangle obtained
by three points $H^{a\bot}_1,H^{a\bot}_2,H^{a\bot}_3$ (resp.
$H^{g\bot}_1,H^{g\bot}_2,H^{g\bot}_3$), and is denoted as
$\triangle(H_1,H_2,H_3)^{a\bot}$ or $\triangle^{a\bot}$ (resp.
$\triangle(H_1,H_2,H_3)^{g\bot}$ or $\triangle^{g\bot}$).
\end{defi}

\begin{rem}\label{3.6} There are two geodesic segments joining $v$
and $w(\ne -v)$ in $\Bbb S^n_H$.  We can choose one geodesic segment
 \textsl{"smaller"} than the other. Here the meaning of  \textsl{"smaller"} is
not  smaller in length (because the length in this model has complex
value) but the one which does not contain two antipodal points.
\end{rem}

The following corollary is an easy consequence of the above definitions.

\begin{cor}\label{3.77} For a triangle in $\Bbb S^2_H$,
\newline (1)
$\triangle(v_1,v_2,v_3)^{a\bot}=\triangle(v^{a\bot}_1,v^{a\bot}_2,v^{a\bot}_3)$.
\newline (2) $(\triangle^{a\bot})^{a\bot}=\triangle$.
\end{cor}

\begin{rem}\label{gg} The relation $(\triangle^{g\bot})^{g\bot}=\triangle$ is not satisfied in general.
In order to get the relation
$(\triangle^{g\bot})^{g\bot}=\triangle$, we have to find a different
type triangle edge and interior construction  for the definition of
$\triangle^{g\bot}$ with the same three vertices.
\end{rem}

Now we can calculate the trigonometric formulas for a triangle in
$\Bbb S^2_H$. We start with any linearly independent non-lightlike
triple ($v_1,v_2,v_3$) of vectors in $\Bbb R^{2,1}$. They determine
a triangle $\triangle(v_1,v_2,v_3)$ formed by {\it smaller}
geodesics. The dual basis of $(v_1,v_2,v_3)$ is another triple
$(w_1,w_2,w_3)$ of vectors in $\Bbb R^{2,1}$, defined by the
conditions $\la v_i,w_i\ra =1$ and $\la v_i,w_j\ra =0$ if $i\ne j$,
for $i,j=1,2,3$. Thus we have
$\triangle(v_1,v_2,v_3)^{a\bot}=\triangle(w_1,w_2,w_3)$. If we let
$V$ and $W$ be the matrices with columns $v_i$ and $w_i$, then they
satisfy the equation $W^t S V=I$, where $S$ is a diagonal matrix
with entries (-1,1,1).  However the matrices of  inner product, $V^t
S V$ and $W^t S W$, are still inverse to each other:
$$ (V^t S V)(W^t S W)= (V^t S V)(V^{-1}  W)= V^t S W=( W^t S V)^t =I.$$

The matrix $ V^t S V$ can be written as

$$  V^t S V =\left(\begin{matrix}

         c_{11}      &c_{12}   &c_{13}  \\
         c_{12}      &c_{22}   &c_{23}  \\
         c_{13}      &c_{23}   &c_{33}  \\
  \end{matrix}\right), \qquad c_{ij}=\la v_i,v_j \ra,
$$
and hence $W^t S W$ is represented as

\begin{equation}\label{66} W^t S W=\frac 1{\det(V^t S V)}   \left(
 \begin{matrix}
         c_{22}c_{33}-c_{23}^2            &c_{13}c_{23}-c_{33}c_{12}     &c_{12}c_{23}-c_{22}c_{13}       \\
         c_{13}c_{23}-c_{33}c_{12}     &c_{11}c_{33}-c_{13}^2            &c_{12}c_{13}-c_{11}c_{23}       \\
         c_{12}c_{23}-c_{22}c_{13}     &c_{12}c_{13}-c_{11}c_{23}    &c_{11}c_{22}-c_{12}^2               \\
  \end{matrix}\right).
\end{equation}

We need another notation. Let's denote the geometric dual of
$\triangle(v_1,v_2,v_3)$ as $\triangle(w'_1,w'_2,w'_3)$.  Then the
formula (\ref{5}) in Corollary \ref{2.11} gives the angle $\theta$
of vertex $v_3$ as
$$ \la w'_1,w'_2\ra = -\vv w'_1\vv \vv w'_2\vv \cos \theta.$$
Here the $w_i$ and $w'_i$ are the same or differ by $-1$, and the difference of
 $\frac {\la w'_1,w'_2\ra}{\vv w'_1\vv \vv w'_2\vv}$ and $\frac {\la w_1,w_2\ra}{\vv w_1\vv \vv w_2\vv}$ is determined by
 $\text{sgn }(\vv w_1 \vv^2\vv w_2\vv ^2)$. Since we get $w_i=w'_i$ if $\vv
 w_i\vv^2<0$, and $w_i=-w'_i$ if $\vv
 w_i\vv^2>0$ from Definition \ref{3.4} or simply
 \begin{equation}\label{x}
 w_i=\text{sgn }(-\vv
 w_i\vv^2)w'_i.
 \end{equation}

 From equation (\ref{66}), we have
$$
\al
\frac {\la w_1,w_2\ra}{\vv w_1\vv \vv w_2\vv}&=\left. \frac{c_{13}c_{23}-c_{33}c_{12}}{\det (V^t S V)}\right/
 \left(  \sqrt{\frac{c_{22}c_{33}-c_{23}^2}{\det (V^t S V)}}
\sqrt{\frac{c_{11}c_{33}-c_{13}^2}{\det (V^t S V)}} \right),\\
&=\text{sgn }\left( (c_{22}c_{33}-c_{23}^2)(c_{11}c_{33}-c_{13}^2)\right) \frac{c_{13}c_{23}-c_{33}c_{12}}{\sqrt{c_{22}c_{33}-c_{23}^2}
\sqrt{c_{11}c_{33}-c_{13}^2}},\\
&=\text{sgn }(\vv w_1\vv^2 \vv w_2\vv^2) \frac{c_{13}c_{23}-c_{33}c_{12}}{\sqrt{c_{22}c_{33}-c_{23}^2} \sqrt{c_{11}c_{33}-c_{13}^2}},
\eal
$$
where we used the fact that $\text{det }(V^t S V)$ is negative and
$\sqrt{\frac{-}{-}}=\frac{\sqrt{-}}{\sqrt{-}}$,
$\sqrt{\frac{+}{-}}=-\frac{\sqrt{+}}{\sqrt{-}}$. Therefore we
conclude
$$
-\cos
\theta=\frac{c_{13}c_{23}-c_{33}c_{12}}{\sqrt{c_{22}c_{33}-c_{23}^2}
\sqrt{c_{11}c_{33}-c_{13}^2}}.
$$

Also we know that $c_{ij}=\la v_i,v_j \ra=\vv v_i\vv \vv v_j\vv
\cosh d_H(v_i,v_j)$ and simply $c_{ij}=\vv v_i\vv \vv v_j\vv \cosh
d_{ij}$  from Theorem \ref{2.10}. So we get
$$
\al
\cos \theta=&\frac{\vv v_1\vv \vv v_2\vv \vv v_3\vv^2 (\cosh d_{12}-\cosh d_{13}\cosh d_{23})}{\sqrt{\vv v_1\vv^2 \vv v_3\vv^2(1-\cosh^2 d_{13})}
 \sqrt{\vv v_2\vv^2 \vv v_3\vv^2(1-\cosh^2 d_{23})}  },\\
=&\frac{\vv v_1\vv \vv v_2\vv \vv v_3\vv^2 (\cosh d_{12}-\cosh d_{13}\cosh d_{23})}{\sqrt{-\vv v_1\vv^2 \vv v_3\vv^2\sinh^2 d_{13}}
                          \sqrt{-\vv v_2\vv^2 \vv v_3\vv^2\sinh^2 d_{23}}  },\\
=&\text{msgn}(-1, \vv v_1\vv^2, \vv v_3\vv^2, \vv v_1\vv^2 \vv w_2\vv^2 \vv v_3\vv^2) \text{msgn}(-1, \vv v_2\vv^2, \vv v_3\vv^2,  \\
 &\vv w_1\vv^2 \vv v_2\vv^2 \vv v_3\vv^2) \frac{\cosh d_{13}\cosh d_{23}-\cosh d_{12}}{\sqrt{\sinh^2 d_{13}}  \sqrt{\sinh^2 d_{23}}  },\\
=&\text{msgn}(-1, \vv v_1\vv^2, \vv v_3\vv^2, \vv v_1\vv^2 \vv w_2\vv^2 \vv v_3\vv^2) \text{msgn}(-1, \vv v_2\vv^2, \vv v_3\vv^2, \\
&\vv w_1\vv^2 \vv v_2\vv^2 \vv v_3\vv^2)  \text{sgn}(\sinh d_{12})
\text{sgn}(\sinh d_{23}) \frac{\cosh d_{13}\cosh d_{23}-\cosh
d_{12}} {\sinh d_{13} \sinh d_{23}}, \eal
$$
by considering $\vv w_1\vv^2 \text{det }(V^t S V)=-\vv v_2\vv^2 \vv
v_3\vv^2\sinh^2 d_{23}$ and $\vv w_2\vv^2$ $\text{det }$ $(V^t S
V)=$ $-\vv v_1\vv^2$ $\vv v_3\vv^2$ $\sinh^2 d_{13}$.

In the above, the function {\it sgn} is defined for pure imaginary
number (for example, sgn($i$)=1 and sgn($-i$)=$-1$), and
$\text{sgn}(\sinh d_{23})$ is negative if and only if $ \vv
v_2\vv^2>0, \vv v_3\vv^2>0,$ and $\vv w_1\vv^2>0$. Then we can show
the following relations by case by case examination.
\begin{equation}\label{7} \al
\text{sgn}(\sinh d_{13})&=-\text{msgn}(-1, -\vv v_1\vv^2, -\vv
v_3\vv^2, -\vv w_2\vv^2)=\text{sgn}(-\vv v_1\vv^3\vv v_3\vv^3\vv
w_2\vv^3) \quad \text{and}\\
\text{sgn}(\sinh d_{23})&=-\text{msgn}(-1, -\vv v_2\vv^2, -\vv
v_3\vv^2, -\vv w_1\vv^2)=\text{sgn}(-\vv v_2\vv^3\vv v_3\vv^3\vv
w_1\vv^3). \eal \end{equation}

 The right hand side of
the equality (\ref{7}) has also negative sign $-1$, when  $ \vv
v_1\vv^2<0, \vv v_3\vv^2<0,$ and $\vv w_2\vv^2<0$.  But we need not
worry about this, because  $ \vv v_1\vv^2<0, \vv v_3\vv^2<0$ implies
$\vv w_2\vv^2>0$. Hence the case does not exist. So the relations
(\ref{7}) are true statements. Therefore we have to simplify the
expression:
$$
\text{msgn}(-1, \vv v_1\vv^2, \vv v_3\vv^2, \vv v_1\vv^2 \vv w_2\vv^2 \vv v_3\vv^2) \text{msgn}
(-1, \vv v_2\vv^2, \vv v_3\vv^2, \vv w_1\vv^2 \vv v_2\vv^2 \vv v_3\vv^2)
$$
$$
\times\ \text{msgn}(-1, -\vv v_1\vv^2, -\vv v_3\vv^2, -\vv w_2\vv^2) \text{msgn}(-1, -\vv v_2\vv^2, -\vv v_3\vv^2, -\vv w_1\vv^2).
$$

\begin{lem}\label{3.8}
$$
\emph{msgn}(-1, \vv v_1\vv^2, \vv v_3\vv^2, \vv v_1\vv^2 \vv
w_2\vv^2 \vv v_3\vv^2) \emph{msgn}(-1, \vv v_2\vv^2, \vv v_3\vv^2,
\vv w_1\vv^2 \vv v_2\vv^2 \vv v_3\vv^2)
$$
$$
\times\ \emph{msgn}(-1, -\vv v_1\vv^2, -\vv v_3\vv^2, -\vv w_2\vv^2)
\emph{msgn}(-1, -\vv v_2\vv^2, -\vv v_3\vv^2, -\vv w_1\vv^2)=1.
$$
\end{lem}

\begin{proof} By using Proposition \ref{3.3}, we see
$$
\al
&\text{msgn}(-1, \vv v_2\vv^2, \vv v_3\vv^2, \vv w_1\vv^2 \vv v_2\vv^2 \vv v_3\vv^2)
\text{msgn}(-1, -\vv v_2\vv^2, -\vv v_3\vv^2, -\vv w_1\vv^2)\\
=&\text{msgn}(-1,-1, \vv v_2\vv^2,-\vv v_2\vv^2, \vv v_3\vv^2,-\vv v_3\vv^2, \vv w_1\vv^2 \vv v_2\vv^2 \vv v_3\vv^2,-\vv w_1\vv^2)\\
&\times\ \text{msgn}(-\vv w_1\vv^2, \vv w_1\vv^2 \vv v_2\vv^2\vv v_3\vv^2)\\
=&\left( \text{msgn}(-\vv w_1\vv^2, \vv w_1\vv^2 \vv v_2\vv^2\vv v_3\vv^2) \right)^2\\
=&1.
\eal
$$
and similarly
$$
\text{msgn}(-1, \vv v_1\vv^2, \vv v_3\vv^2, \vv v_1\vv^2 \vv w_2\vv^2 \vv v_3\vv^2)
\text{msgn}(-1, -\vv v_1\vv^2, -\vv v_3\vv^2, -\vv w_2\vv^2)=1.
$$
\end{proof}

As is shown, we conclude
$$
\cos \theta=\frac{\cosh d_{13}\cosh d_{23}-\cosh d_{12}}{\sinh d_{13} \sinh d_{23}}.
$$

Letting $A,B,C$ stand for the angles at $v_1,v_2,v_3$ and $a,b,c$
for the extended hyperbolic lengths of opposite sides,  we obtain
the hyperbolic law of cosine on the hyperbolic sphere $\Bbb S^2_H$:
$$
\cos C=\frac{\cosh a\cosh b-\cosh c}{\sinh a \sinh b}.
$$
Also we can easily deduce the spherical law of cosine on the
spherical sphere $\Bbb S^2_S$ by using $i\cdot d_S=d_H$ and so
$\cosh d_H=\cos d_S,$ $\sinh d_H=i\sin d_S$, where $a,b,c$ represent
the extended spherical length,
$$
\al
\cos C&=\frac{\cosh (ai)\cosh (bi)-\cosh (ci)}{\sinh (ai)\sinh (bi)},\\
      &=\frac{\cos c-\cos a\cos b}{\sin a \sin b}.
\eal
$$

To obtain the dual cosine law, we start our argument from a triangle
$\triangle(v_1,v_2,v_3) $ with its geometric dual
$\triangle^{g\bot}$ written by $\triangle(w_1,w_2,w_3)$.  In the
proof of cosine law, $\triangle(w_1,w_2,w_3)$ means an algebraic
dual, but from now $\triangle(w_1,w_2,w_3)$ denotes a geometric dual
for convenience. The angles and edges of $\triangle$ and
$\triangle^{g\bot}$ are shown in Fig. 10.

Lemma \ref{2.4} and Definition \ref{2.2} deduce the relations $-\cos
A=\cosh a',$ $ -\cos B=\cosh b',$ and $-\cos C=\cosh c'$, but do not
gives the relations $-\cos A'=\cosh a, -\cos B'=\cosh b,$ and $-\cos
C'=\cosh c$ by Remark \ref{gg}. By comparison of $\triangle^{a\bot}$
and $\triangle^{g\bot}$, and comparison of
$(\triangle^{a\bot})^{g\bot}$ and
$(\triangle^{a\bot})^{a\bot}=\triangle$, and the relation (\ref{x}),
we can get the exact relations between $\cos A', \cos B', \cos C'$
and $\cosh a, \cosh b, \cosh c$:
$$
\al
-\cos A'=&\cosh a\ \text{sgn}(\vv v_2\vv^2 \vv v_3\vv^2 \vv w_2\vv^2 \vv w_3\vv^2),\\
-\cos B'=&\cosh b\ \text{sgn}(\vv v_1\vv^2 \vv v_3\vv^2 \vv w_1\vv^2 \vv w_3\vv^2),\\
-\cos C'=&\cosh c\ \text{sgn}(\vv v_1\vv^2 \vv v_2\vv^2 \vv w_1\vv^2 \vv w_2\vv^2).
\eal
$$

\begin{figure}[h]
\begin{center}
\includegraphics[width=0.55\textwidth]{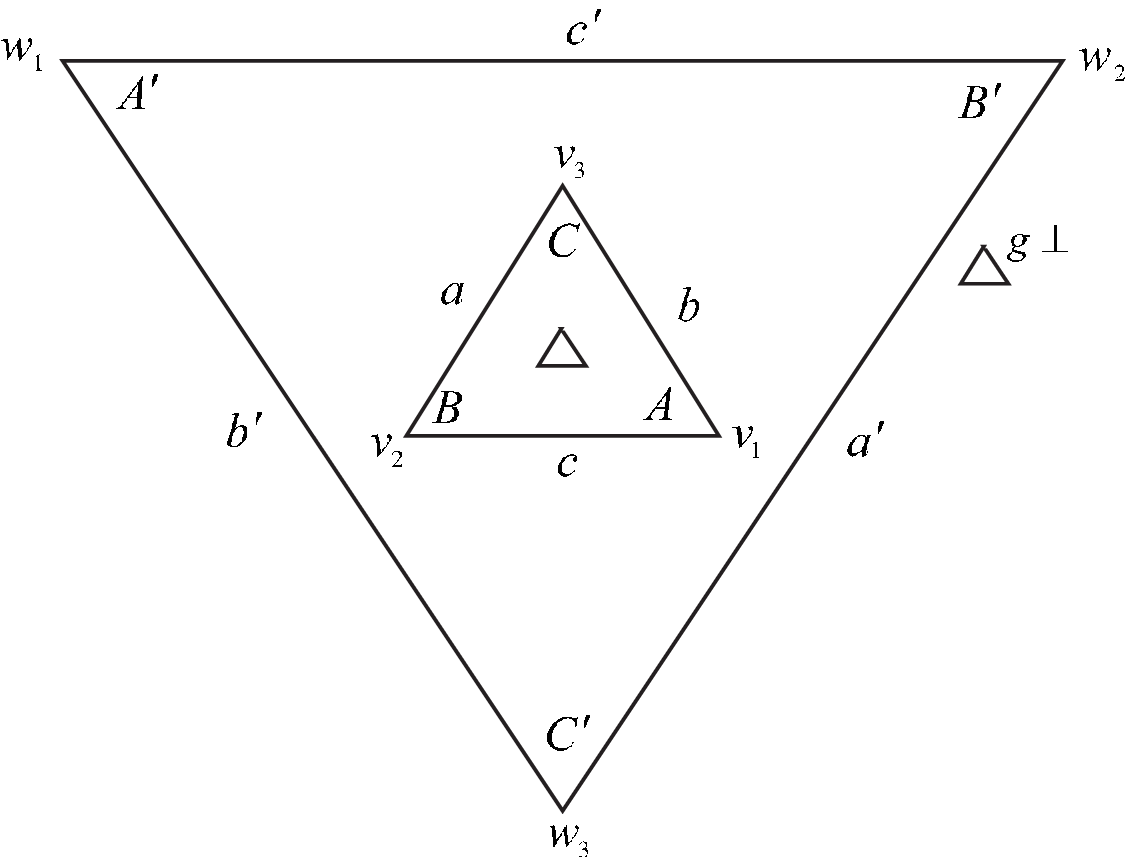}
\caption{\textbf{9}}
\end{center}
\end{figure}

We already get the cosine law which is adapted to the triangle
$\triangle^{g\bot}$:
\begin{equation}\label{8} \cos C'=\frac{\cosh
a'\cosh b'-\cosh c'}{\sinh a' \sinh b'}. \end{equation}
 The
formula (\ref{8}) is changed to
$$
-\cosh c\ \text{sgn}(\vv v_1\vv^2 \vv v_2\vv^2 \vv w_1\vv^2 \vv w_2\vv^2)=\frac{\cos A\cos B+\cos C}{\sinh a' \sinh b'},
$$
and we have to show that
\begin{equation}\label{9}
\sin A \sin B=-\sinh a' \sinh b' \ \text{sgn}(\vv v_1\vv^2 \vv
v_2\vv^2 \vv w_1\vv^2 \vv w_2\vv^2)
\end{equation}
in order to obtain the dual cosine law,
$$
\cosh c=\frac{\cos A\cos B+\cos C}{\sin A \sin B}.
$$

The above formula (\ref{9}) is also changed to
$$
\al
&\sin A \sin B\ \text{sgn}(\vv v_1\vv^2 \vv v_2\vv^2 \vv w_1\vv^2 \vv w_2\vv^2)\\
=&-\text{sgn}(\sinh a') \text{sgn}(\sinh b')\sqrt{\sinh^2 a'}\sqrt{\sinh^2 b'}\\
=&-\text{sgn}(\sinh a') \text{sgn}(\sinh b')\sqrt{-1+\cos^2 A}\sqrt{-1+\cos^2 B}\\
=&\text{sgn}(\sinh a') \text{sgn}(\sinh b')\text{msgn}(-1,\sin^2 A) \text{msgn}(-1,\sin^2 B)\sqrt{\sin^2 A}\sqrt{\sin^2 B},
\eal
$$
hence we need
$$
\al
&\text{sgn}(\vv v_1\vv^2 \vv v_2\vv^2 \vv w_1\vv^2 \vv w_2\vv^2)\\
=&\text{sgn}(\sinh a') \text{sgn}(\sinh b')\text{msgn}(-1,\sin^2 A) \text{msgn}(-1,\sin^2 B)\text{sgn}(\sin A)\text{sgn}(\sin B).
\eal
$$

For complex numbers $z_1$ and $z_2$, if there exists a positive
number $\alpha$ such that $z_1=\alpha z_2$, then let's denote simply
as $z_1\sim z_2$. Then by easy checking, we know $\sinh a'\sim -\vv
v_1\vv^3\vv w_2\vv^3\vv w_3\vv^3$ and $\sin A\sim -i\vv v_1\vv \vv
w_2\vv \vv w_3\vv$. Also we can easily find $$\text{sgn}(-\vv
v_1\vv^3\vv w_2\vv^3\vv w_3\vv^3)=\text{sgn}(-i\vv v_1\vv \vv w_2\vv
\vv w_3\vv)=-\text{msgn}(-1,-\vv v_1\vv^2, -\vv w_2\vv^2, -\vv
w_3\vv^2).$$

Therefore we can get the following identities:
$$
\al
\text{sgn}(\sinh a')&=\text{sgn}(\sin A)=-\text{msgn}(-1,-\vv v_1\vv^2, -\vv w_2\vv^2, -\vv w_3\vv^2),\\
\text{sgn}(\sinh b')&=\text{sgn}(\sin B)=-\text{msgn}(-1,-\vv v_2\vv^2, -\vv w_1\vv^2, -\vv w_3\vv^2),
\eal
$$
and
$$
\al
\text{sgn}(\sin^2 A)&=\text{sgn}(-\vv v_1\vv^2\vv w_2\vv^2\vv w_3\vv^2),\\
\text{sgn}(\sin^2 B)&=\text{sgn}(-\vv v_2\vv^2\vv w_1\vv^2\vv
w_3\vv^2). \eal
$$

The only thing left to show is the following lemma.

\begin{lem}\label{3.99} $\emph{msgn}(-1,-\vv v_1\vv^2\vv
w_2\vv^2\vv w_3\vv^2) \emph{msgn}(-1,-\vv v_2\vv^2\vv w_1\vv^2\vv
w_3\vv^2)$ \newline $\text{       }=\emph{sgn}(\vv v_1\vv^2 \vv
v_2\vv^2 \vv w_1\vv^2 \vv w_2\vv^2). $
\end{lem}

\begin{proof} The left  hand side of the above equality is equal
to
$$
\al
&\text{msgn}(-1,-1,-\vv v_1\vv^2\vv w_2\vv^2\vv w_3\vv^2,-\vv v_2\vv^2\vv w_1\vv^2\vv w_3\vv^2)\\
&\times\ \text{msgn}(\vv v_1\vv^2\vv w_2\vv^2\vv w_3\vv^2,\vv v_2\vv^2 \vv w_1\vv^2 \vv w_3\vv^2)\\
=&\text{msgn}(-1,-1,-\vv v_1\vv^2\vv w_2\vv^2\vv w_3\vv^2,-\vv v_2\vv^2\vv w_1\vv^2\vv w_3\vv^2,\vv v_1\vv^2\vv w_2\vv^2\vv w_3\vv^2\\
&,\vv v_2\vv^2 \vv w_1\vv^2 \vv w_3\vv^2)\times\
\text{msgn}(\vv v_1\vv^2 \vv v_2\vv^2 \vv w_1\vv^2 \vv w_2\vv^2,\vv v_1\vv^2 \vv v_2\vv^2 \vv w_1\vv^2 \vv w_2\vv^2)\\
=&\text{sgn}(\vv v_1\vv^2 \vv v_2\vv^2 \vv w_1\vv^2 \vv w_2\vv^2),
\eal
$$
where we used Proposition \ref{3.3} b),c), and d).
\end{proof}

As is shown, we deduce the dual cosine law on the hyperbolic sphere $\Bbb S^2_H$.

In order to get the dual cosine law on the spherical sphere $\Bbb
S^2_S$, we need only $i\cdot d_S=d_H$ as before.

We considered only triangle with non-lightlike vertex vectors, i.e.,
without ideal vertices. If we permit lightlike vector, then the
values $\vv \cdot\vv$ become $0$ and angles and lengths can be $0$
or $\infty$. Even in this degenerated case, we can easily convince
the cosine and dual cosine law, too. Therefore we can summarize the
cosine law and dual cosine law for $\Bbb S^2_H$ and $\Bbb S^2_S$ in
the following theorem.

\begin{thm}\label{3.10} Letting $A,B,C$ stand for the angles and
$a,b,c$ for the extended hyperbolic  lengths of opposite sides of a
given triangle, we obtain the hyperbolic cosine law and the dual
cosine law on the hyperbolic sphere $\Bbb S^2_H$,
$$
\al
\cos C&=\frac{\cosh a\cosh b-\cosh c}{\sinh a \sinh b},\\
\cosh c&=\frac{\cos A\cos B+\cos C}{\sin A \sin B}.
\eal
$$
Also we have the spherical cosine law and dual cosine law on the
spherical sphere $\Bbb S^2_S$, where $a,b,c$ represent the extended
spherical lengths,
$$
\al
\cos C&=\frac{\cos c-\cos a\cos b}{\sin a \sin b},\\
\cos c&=\frac{\cos A\cos B+\cos C}{\sin A \sin B}.
\eal
$$
\end{thm}

Now we consider the cosine laws on the hyperbolic sphere $\Bbb
S^3_H$ or the spherical sphere $\Bbb S^3_S$, then we should consider
two more cases of triangles.

If a hyperplane containing the triangle does not intersect to
$\partial\Bbb H^3$, then we can send this triangle to the equator
(=$\Bbb S^3_H\cap\{x|x_0=0\}$) of $\Bbb S^3_H$ by an isometry. Hence
the distance becomes $i$ times the distance on the standard
Euclidean unit sphere. Therefore the above result of Theorem
\ref{3.10} also satisfied by the well known spherical trigonometry.

\begin{figure}[h]
\begin{center}
\includegraphics[width=0.8\textwidth]{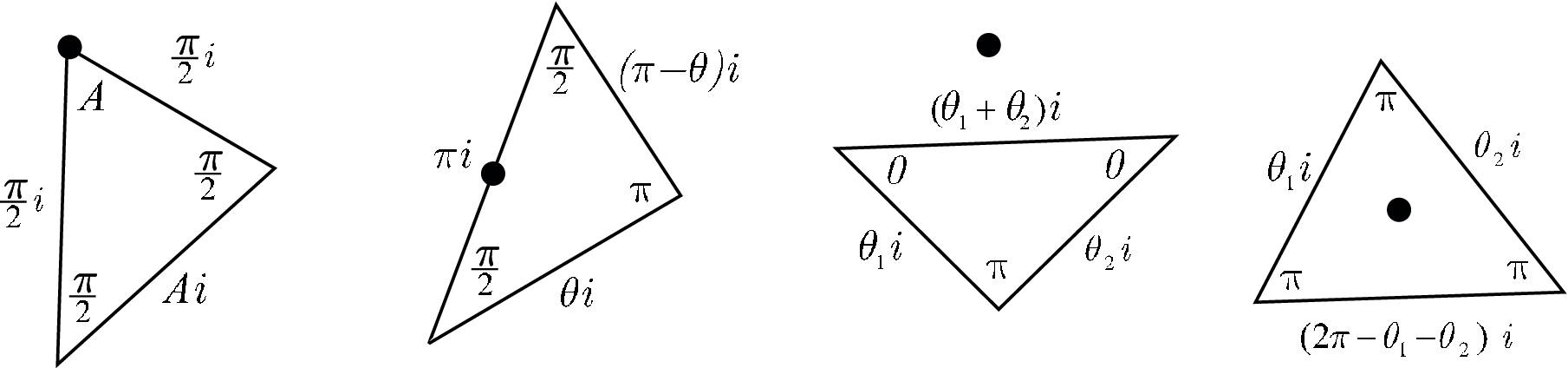}
\caption{\textbf{11}}
\end{center}
\end{figure}

If a hyperplane containing the triangle is tangent to $\partial\Bbb
H^3$, then there are only four types of triangles in the sense of
the intersection of three hemispheres (see Fig. 11).

These four types also satisfy the above cosine laws in some sense,
if we permit $\frac{0}{0}=$ a certain number or move the denominator
of the formulas to the other side.

Therefore we can conclude a theorem about all kinds of triangles in
$\Bbb S^n_H$ or $\Bbb S^n_S$.

\begin{thm}\label{3.101} For a given triangle in $\Bbb S^n_H$ (resp. $\Bbb
S^n_S$), the triangle satisfies the hyperbolic (resp. spherical)
 cosine and dual cosine laws as in Theorem \ref{3.10}.
\end{thm}

\vtwo \leftline{\bf 4.2 Sine law} \vtwo

The hyperbolic sine law is easily obtained by the following steps.
First we assume that all vertices of a triangle are not ideal
vertices. From the dual cosine law for a right triangle with
$C=\frac{\pi}2$, we have
\begin{equation}\label{10}
\cosh b=\frac{\cos A\cos C+\cos B}{\sin A \sin C}=\frac{\cos B}{\sin
A},
\end{equation}
and also cosine law induces $\cosh c=\cosh a\cosh b$ and
\begin{equation}\label{11}
\cos B=\frac{\cosh a\cosh c-\cosh b}{\sinh a \sinh c}.
\end{equation}
By substituting (\ref{10}) and $\cosh c=\cosh a\cosh b$ into
(\ref{11}), we get
$$
\sinh a= \sin A\sinh c.
$$

Now given any triangle with sides $\tilde a,\tilde b,\tilde c$ and
angles $\tilde A,\tilde B,\tilde C$, the altitude $h$ corresponding
to side $a$, so we can induce $\sinh h=\sin \tilde A \sinh \tilde b$
and also $\sinh h=\sin \tilde B \sinh \tilde a$. Here altitude line
can be constructed by joining one vertex point and the dual point of
the line which passes the other two points. In the proof, the
non-ideal vertex condition is necessary used for cancelation. When
we consider ideal vertex case, then the sine law also satisfied by
easy checking. \noindent Therefore we proved the following theorem
for hyperbolic sine law and spherical sine law.

\begin{thm}\label{3.11} Letting $A,B,C$ stand for the angles and
$a,b,c$ for the extended hyperbolic lengths of opposite sides of a
given triangle, we obtain the hyperbolic sine law on the hyperbolic
sphere $\Bbb S^2_H$,
$$
\frac{\sinh a}{\sin A}=\frac{\sinh b}{\sin B}=\frac{\sinh c}{\sin C}.
$$
Also we have the spherical sine law on the spherical sphere $\Bbb
S^2_S$, where $a,b,c$ represent the extended spherical lengths,
$$
\frac{\sin a}{\sin A}=\frac{\sin b}{\sin B}=\frac{\sin c}{\sin C}.
$$
\end{thm}
We introduce another proof.

\begin{proof} We know
$$\frac{\sinh^2 a}{\sin^2 A}=\frac{\sinh^2 a \sinh^2 b \sinh^2
c}{1-\cosh^2 a -\cosh^2 b -\cosh^2 c+2\cosh a \cosh b \cosh c},$$ so
we conclude $\frac{\sinh^2 a}{\sin^2 A}=\frac{\sinh^2 b}{\sin^2
B}=\frac{\sinh^2 c}{\sin^2 C},$ in particular, it is also satisfied
when a denominator or numerator of the formula takes $0$ or
$\infty$. Now we have to show that  $\frac{\sinh a}{\sin
A}=\frac{\sinh b}{\sin B}=\frac{\sinh c}{\sin C}$ for a non-ideal
vertices triangle.

From $\sinh^2 a \sin^2 B=\sinh^2 b \sin^2 A$, it follows
continuously that $$\sqrt{\sinh^2 a \sin^2 B}=\sqrt{\sinh^2 b \sin^2
A}$$
$$\text{msgn}(\sinh^2 a, \sin^2 B)\sqrt{\sinh^2 a}\sqrt{\sin^2
B}=\text{msgn}(\sinh^2 b, \sin^2 A)\sqrt{\sinh^2 b}\sqrt{\sin^2 A}$$
$$\text{msgn}(\sinh^2 a, \sin^2 B)\text{sgn}(\sinh a)\text{sgn}(\sin
B)\sinh a \sin B$$ $$=\text{msgn}(\sinh^2 b, \sin^2
A)\text{sgn}(\sinh b)\text{sgn}(\sin A)\sinh b \sin A.$$ Hence if

\noindent$$\text{msgn}(\sinh^2 a, \sin^2 B)\text{sgn}(\sinh
a)\text{sgn}(\sin B)=\text{msgn}(\sinh^2 b, \sin^2
A)\text{sgn}(\sinh b)\text{sgn}(\sin A)$$ is satisfied, then the
proof ends.

We already know that

\noindent$\text{sgn}(\sinh a)=-\text{msgn}(-1,-\vv v_2\vv^2, -\vv
v_3\vv^2, -\vv w_1\vv^2)$, $\text{sgn}(\sinh^2 a)=\text{sgn}(\vv
v_2\vv^2\vv v_3\vv^2\vv w_1\vv^2)$,

\noindent$\text{sgn}(\sin A)=-\text{msgn}(-1,-\vv v_1\vv^2, -\vv
w_2\vv^2, -\vv w_3\vv^2)$, $\text{sgn}(\sin^2 A)=\text{sgn}(-\vv
v_1\vv^2\vv w_2\vv^2\vv w_3\vv^2).$

\noindent It suffices to show that the following formula is true.

$$\text{msgn}(-1,-\vv v_1\vv^2, -\vv w_2\vv^2, -\vv
w_3\vv^2)\text{msgn}(-1,-\vv v_1\vv^2, -\vv v_3\vv^2, -\vv
w_2\vv^2)$$ $$\times \text{msgn}(-1,-\vv v_2\vv^2, -\vv w_1\vv^2,
-\vv w_3\vv^2)\text{msgn}(-1,-\vv v_2\vv^2, -\vv v_3\vv^2, -\vv
w_1\vv^2)$$ $$\times \text{msgn}(-\vv v_1\vv^2\vv w_2\vv^2\vv
w_3\vv^2, \vv v_1\vv^2\vv v_3\vv^2\vv w_2\vv^2)\text{msgn}(-\vv
v_2\vv^2 \vv w_1\vv^2 \vv w_3\vv^2, \vv v_2\vv^2\vv v_3\vv^2\vv
w_1\vv^2)=1$$

We left the proof of the $\emph{msgn}$ equality as an easy exercise
for readers.
\end{proof}

As is shown in Theorem \ref{3.101}, we also similarly induce the
sine law on $\Bbb S^n_H$ or $\Bbb S^n_S$ (easy check). Therefore we
can conclude the following theorem.

\begin{thm} For a given triangle in $\Bbb S^n_H$ (resp. $\Bbb
S^n_S$), the triangle satisfies the hyperbolic (resp. spherical)
 sine law as in Theorem \ref{3.11}.
\end{thm}

\vtwo \leftline{\bf 4.3 Applications for hyperbolic polygons} \vtwo

\begin{figure}[h]
\begin{center}
\includegraphics[width=0.6\textwidth]{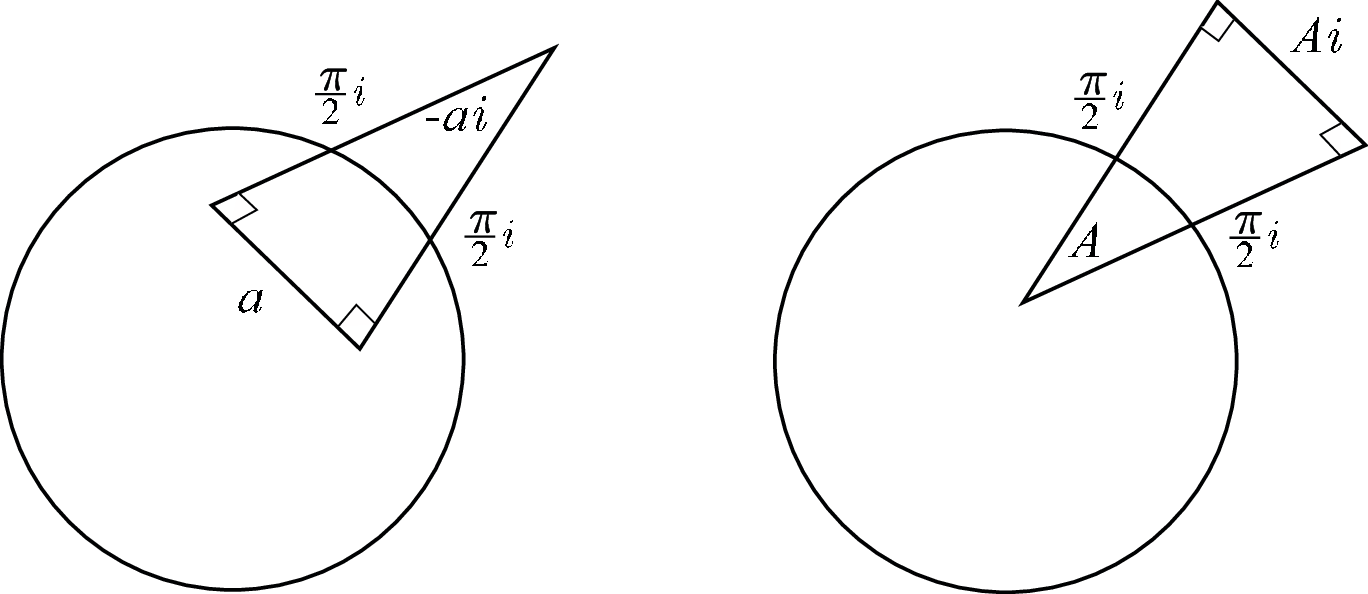}
\caption{\textbf{12}}
\end{center}
\end{figure}

There  are many formulas for Lambert quadrilaterals and pentagons
and hexagons, and these were  shown by an unified  method which
starts from a rectangular hexagon in Fenchel's book \cite{4}. Also
the above general version of cosine laws and sine law also induce
those formulas about all polygons which was mentioned in \cite{4}.
Readers can notice that our interpretation gives an easy and natural
way to understand.

Here we only need the definition of angle and the fact that the
distance between $x$ and $x^{\bot}$ is $\frac{\pi}2 i$ (see Fig.
12). Also it is convenient to remember that
$$\sin ix=i \sinh x,~~ \sinh ix=i \sin x,~~ \sinh (x+\pi i)=- \sinh x,~~ \sinh (x+\frac{\pi}{2}i)=i \cosh x,$$
$$\cos ix=\cosh x,~~ \cosh ix= \cos x,~~ \cosh (x+\pi i)=- \cosh x,~~ \cosh (x+\frac{\pi}{2}i)=i \sinh x.$$

Now we examine four special cases. At first Consider a quadrilateral
with consecutive two right angles shown in Fig. 13 below. We know
that the lengths between 1,2 and 2,3 and 3,1 are $a+\frac{\pi}2 i$,
$c$ and $b+\frac{\pi}2 i$ respectively, and the angles at 1, 2 and 3
are $-d i$, $B$ and $A$.

\begin{figure}[h]
\begin{center}
\includegraphics[width=0.4\textwidth]{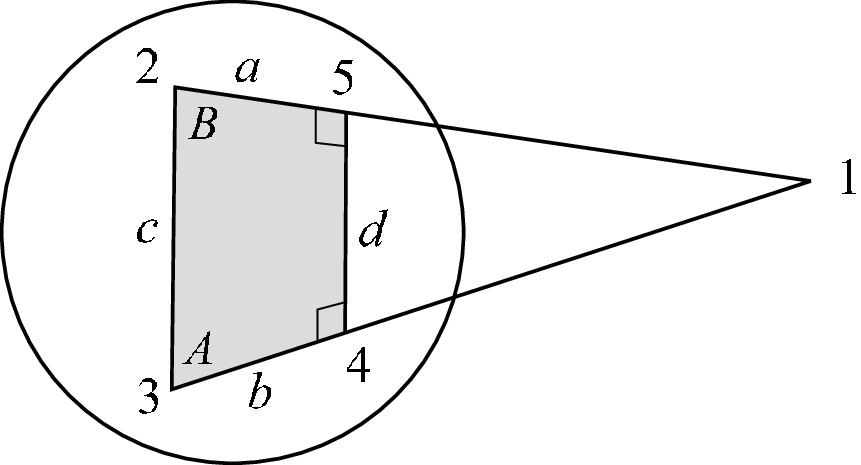}
\caption{\textbf{13}}
\end{center}
\end{figure}

 \noindent From the cosine law for a triangle
(1,2,3),  we obtain formulas:
$$
\al
\cos (-di)&=\frac{\cosh (a+\frac{\pi}2 i)\cosh (b+\frac{\pi}2 i)-\cosh c}{\sinh (a+\frac{\pi}2 i) \sinh (b+\frac{\pi}2 i)}
 \quad \rightarrow\quad \cosh d&=\frac{\sinh a\sinh b+\cosh c}{\cosh a \cosh b},\\
\cos A&=\frac{\cosh c\cosh (b+\frac{\pi}2 i)-\cosh (a+\frac{\pi}2
i)}{\sinh c\sinh (b+\frac{\pi}2 i)}
 \quad \rightarrow\quad\cos A&=\frac{\cosh c\sinh b-\sinh a}{\sinh c \cosh b}.\\
\eal
$$
Also the dual cosine law induces
$$
\al \cosh (a+\frac{\pi}2 i)&=\frac{\cos B\cos (-di)+\cos A}{\sin B
\sin (-di)}
 \quad \rightarrow\quad \sinh a&=\frac{\cos B\cosh d+\cos A}{\sin B
\sinh d},\\
\cosh c&=\frac{\cos A\cos B+\cos (-di)}{\sin A\sin B}
 \quad \rightarrow\quad\cosh c&=\frac{\cos A\cos B+\cosh d}{\sin A\sin B}.\\
\eal
$$

And the sine law implies
$$\frac{\sinh (a+\frac{\pi}2 i)}{\sin A}=\frac{\sinh (b+\frac{\pi}2 i)}{\sin B}=\frac{\sinh c}{\sin (-di)}\quad \rightarrow\quad
  \frac{\cosh a}{\sin A}=\frac{\cosh b}{\sin B}=\frac{\sinh c}{\sinh d}.$$

\begin{figure}[h]
\begin{center}
\includegraphics[width=0.3\textwidth]{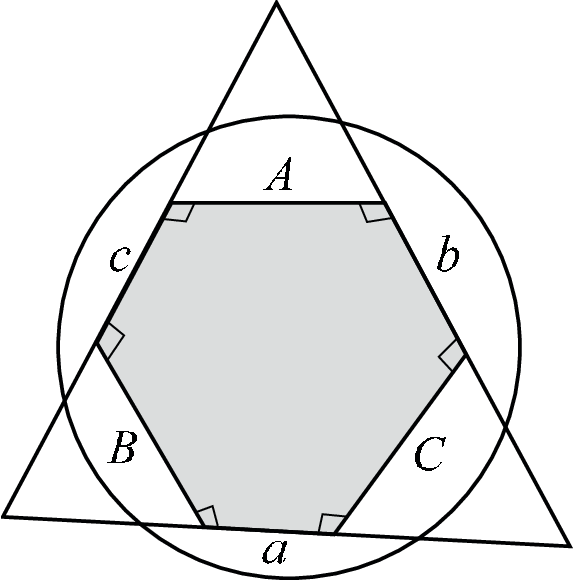}
\caption{\textbf{14}}
\end{center}
\end{figure}

A rectangular hyperbolic hexagon can be perceived as a truncated
triangle (see Fig. 14). So the triangle has  lengths $a+\pi i$,
$b+\pi i$ and $c+\pi i$, and angles $-A i$,$-B i$ and $-C i$.

\noindent From the cosine law and dual cosine law, we get
$$
\cosh C=\frac{\cosh a\cosh b+\cosh c}{\sinh a \sinh b}
$$
and
$$
\cosh c=\frac{\cosh A\cosh B+\cosh C}{\sinh A \sinh B}.
$$

\noindent The sine law shows
$$\frac{\sinh a}{\sinh A}=\frac{\sinh b}{\sinh B}=\frac{\sinh c}{\sinh C}.$$

\begin{figure}[h]
\begin{center}
\includegraphics[width=0.35\textwidth]{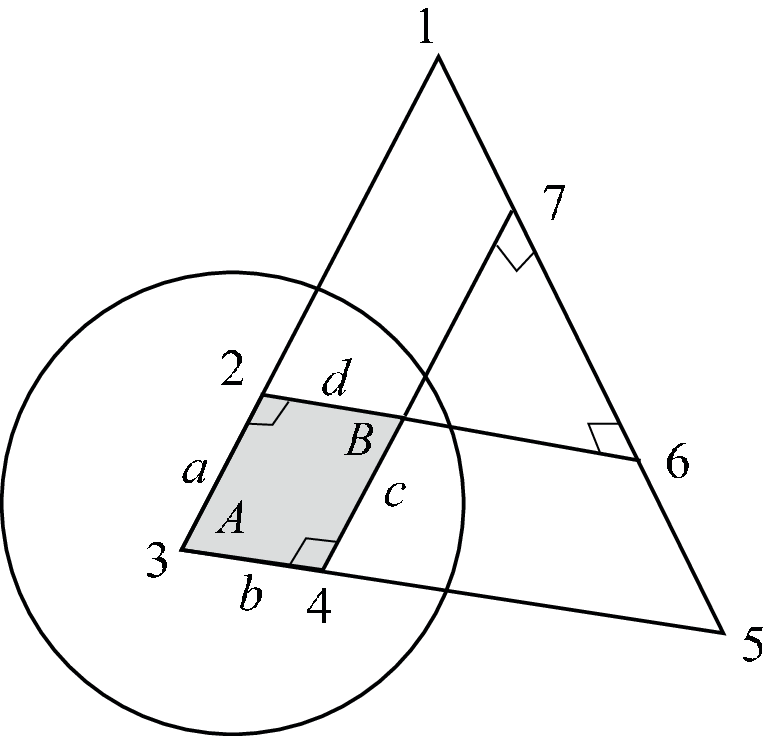}
\caption{\textbf{15}}
\end{center}
\end{figure}

A hyperbolic quadrilateral with two opposite rectangular angles also
can be applicable (see Fig. 15). We know that
$$\overline{13}=a+\frac{\pi}2 i, \quad\overline{35}=b+\frac{\pi}2 i,
\quad\overline{15}=\overline{16}+\overline{57}-\overline{67}=\frac{\pi}2
i+\frac{\pi}2 i-B i,$$ and $$\angle 1=-i\cdot
\overline{26}=-i(d+\frac{\pi}2 i)=\frac{\pi}2 -d i, \quad\angle 5
=\frac{\pi}2 -c i, \quad\angle 3=A.$$
\noindent Hence by generalized
hyperbolic cosine and sine laws, we get
$$
\cos A=\frac{\sinh a\sinh b-\cos B}{\cosh a \cosh b}.
$$
$$
\sinh a=\frac{\cos A\sinh d+\sinh c}{\sin A \cosh d},
$$
$$\frac{\sin B}{\sin A}=\frac{\cosh a}{\cosh c}=\frac{\cosh b}{\cosh d}.$$

\begin{figure}[h]
\begin{center}
\includegraphics[width=0.35\textwidth]{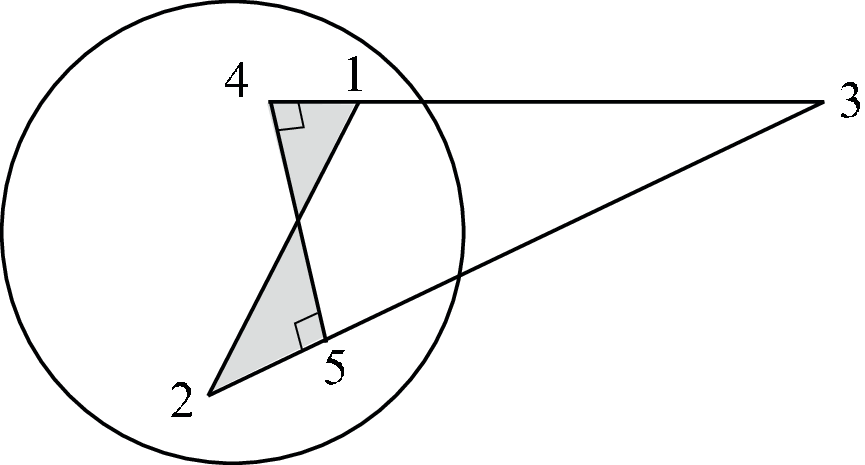}
\caption{\textbf{16}}
\end{center}
\end{figure}

Even in a self intersecting quadrilateral (see Fig. 16), we can
apply the generalized hyperbolic trigonometry. From the general
triangle $\triangle(1,2,3)$, we will get the trigonometry of the
quadrilateral $(1,2,5,4).$

One can easily examine the other formulas for various hyperbolic
polygons in the similar way.

\vtwo \leftline{\bf 4.4 Applications for de Sitter polygons} \vtwo

\begin{figure}[h]
\begin{center}
\includegraphics[width=0.6\textwidth]{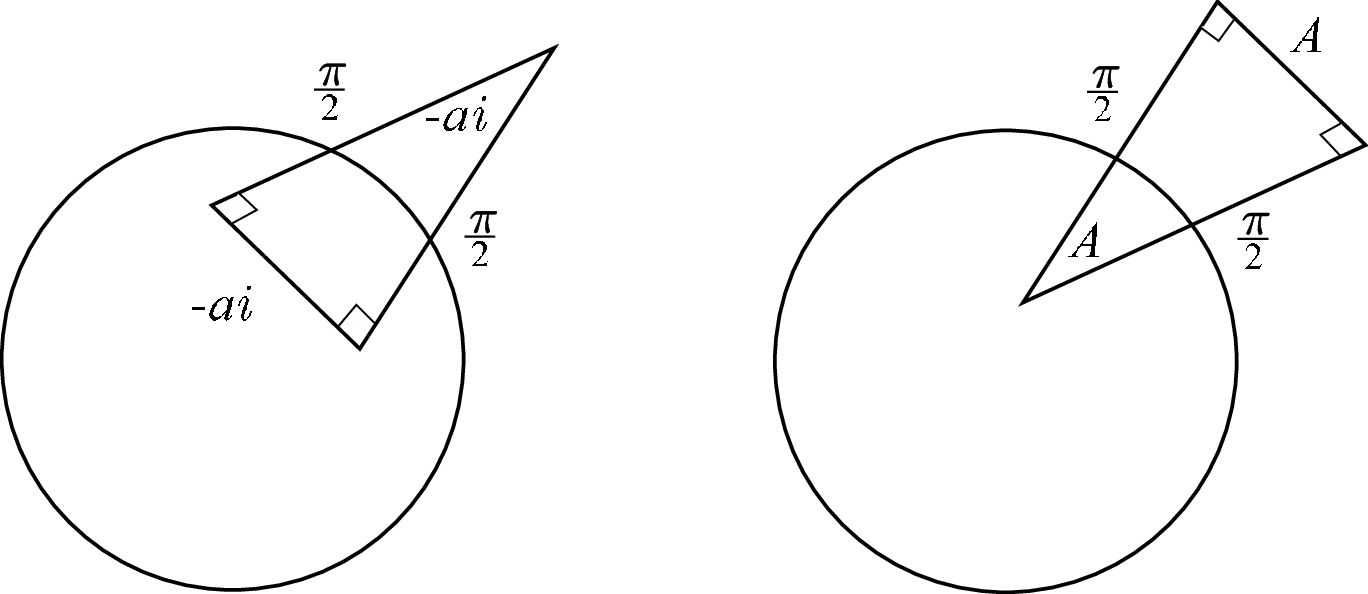}
\caption{\textbf{17}}
\end{center}
\end{figure}

The generalized spherical cosine and sine laws can be used for the
polygons on $S^n_1$. So we  can get many formulas for polygons on
$S^n_1$ by the similar way of \S 4.3, those formulas are not unknown
yet as I know. Especially for a triangle contained in $S^n_1$, Dzan
\cite{32} also induced the same spherical type cosine and sine laws.

First, we need the basic facts about lengths and angles: The
distance between $x$ and $x^{\bot}$ is $\frac{\pi}2$, and the angle
$\angle (x_p,y_p)$ is $d_S(x,y)$ (see Fig. 17). We have to define a
timelike (resp. spacelike) edge as the geodesic edge whose tangent
vector is timelike (resp. spacelike) vector, then we know that a
time edge inside of the Lorentzian part has positive pure imaginary
length and a space edge has positive real length on the extended de
Sitter space (see Convention \ref{1.4}).

\begin{figure}[h]
\begin{center}
\includegraphics[width=0.4\textwidth]{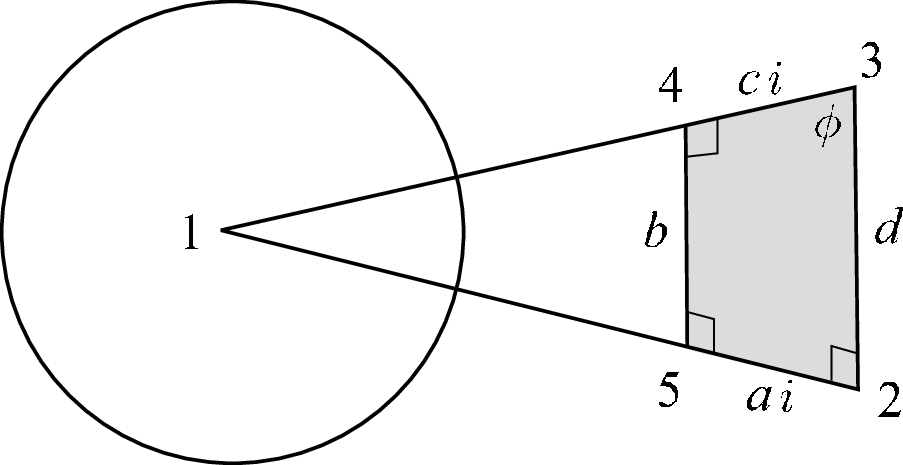}
\caption{\textbf{18}}
\end{center}
\end{figure}

Now we examine the two cases. First one is a Lambert quadrilateral
shown in  Fig. 18.  We know that  the lengths between 1,2 and 1,3
and 2,3 are $ai+\frac{\pi}2$ and $ci+\frac{\pi}2$ and $d$
respectively, and the angles at 1,2 and 3 are $b,\frac {\pi}2$ and
$\phi$. Here $\phi$ is a complex number and all the others are
positive real numbers. From the spherical cosine law for a right
triangle (1,2,3), we obtain three formulas
\begin{equation}\label{12}
\aligned
\cos b&=\frac{\sinh a\sinh c+\cos d}{\cosh a \cosh c},\\
\cos \phi&=-i \frac{\sinh a-\sinh c\cos d}{\cosh c\sin d},\\
\sinh c&=\cos d\sinh a.
\endaligned
\end{equation}
 Also the spherical dual cosine law induces
\begin{equation}\label{13} \cos b=\cos d\sin \phi,\   \cos\phi=-i
\sinh a\sin b, \text{  and  } \sinh a=i \cot b \cot
\phi.\end{equation} The formulas (\ref{12}) and the middle one of
(\ref{13}) induce the inequality $\sinh a>\sinh c\cos d$. And the
sine law gives us
$$\frac{\sin d}{\sin b}=\cosh c=\frac{\cosh a}{\sin \phi}.$$

A pentagon with four right angles in the de Sitter space can be
perceived as a truncated triangle (see Fig. 19).  From the figure,
the triangle (1,2,3) has three side of lengths $ai+\frac{\pi}2 $,
$ei+\frac{\pi}2 $ and $ci+\pi$ and three angles $\phi,b$ and $d$.
Here $\phi$ is a complex number and all the others are positive real
numbers. So we get six formulas from the spherical cosine and dual
cosine laws.

\begin{figure}[h]
\begin{center}
\includegraphics[width=0.5\textwidth]{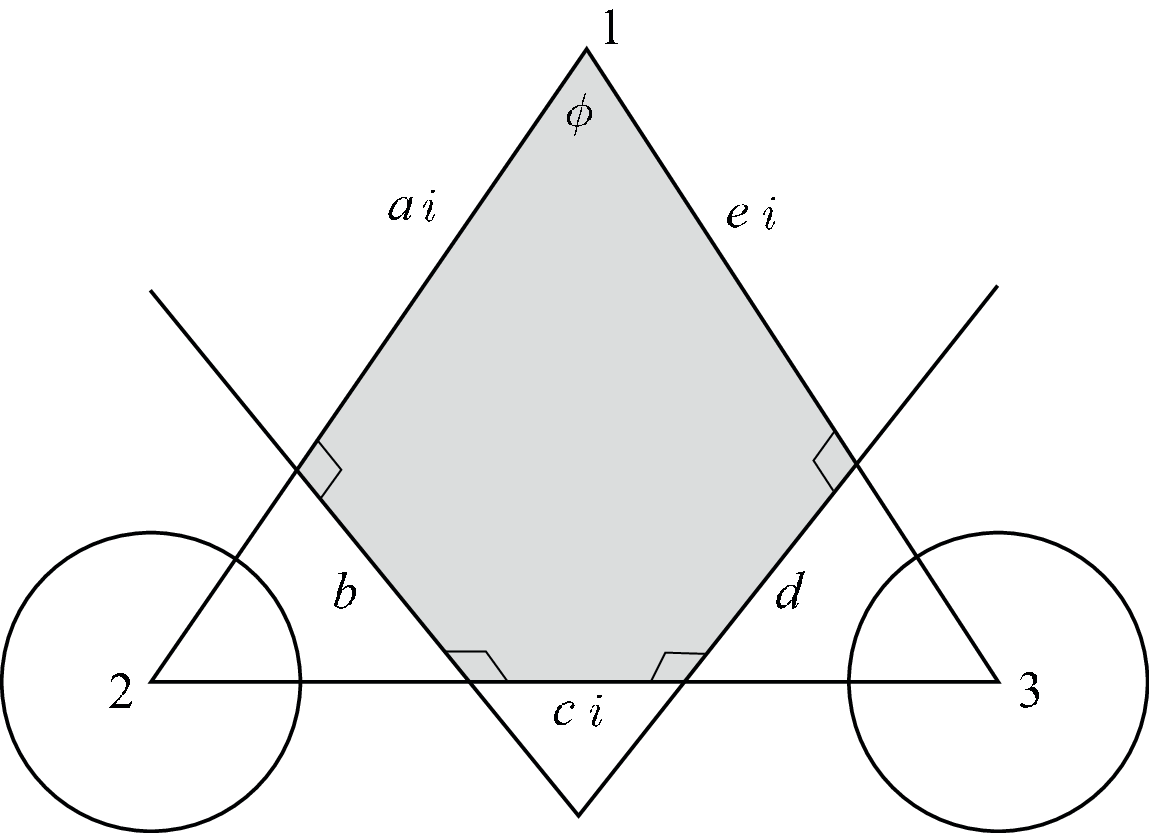}
\caption{\textbf{19}}
\end{center}
\end{figure}

They are
\begin{equation}\label{14}
\aligned
\cos b&=\frac{\sinh a\cosh c+\sinh e}{\cosh a \sinh c},\\
\cos d&=\frac{\sinh e\cosh c+\sinh a}{\cosh e \sinh c},\\
\cos \phi&= \frac{\sinh a\sinh e-\cosh c}{\cosh a\cosh e},
\endaligned
\end{equation}
 and
\begin{equation}\label{15}
\aligned
-\cosh c&=\frac{\cos b\cos d+\cos \phi}{\sin b \sin d},\\
-i \sinh a&=\frac{\cos b\cos \phi+\cos d}{\sin b \sin\phi},\\
-i \sinh e&=\frac{\cos d\cos \phi+\cos b}{\sin d \sin\phi}.
\endaligned
\end{equation}
We can easily show that the angle $b$ and $d$ are smaller than
$\frac{\pi}{2}$, so we get $\cos b, \cos d>0$. Then from the first
formula of (\ref{15}), we have $\cos \phi<0$, so the third formula
of (\ref{14}) gives us an inequality, $\sinh a\sinh e<\cosh c$. And
the sine law implies
$$-i \frac{\sinh c}{\sin \phi}=\frac{\cosh e}{\sin b}=\frac{\cosh a}{\sin d}.$$

Readers can easily induce the trigonometry formula for some de
Sitter polygons with six variables of lengths and angles and the
others rectangular angles by the similar way.

\vskip 1pc
 Lastly we want to remark some problems. Even though the
properties on the extended space are very natural, our proof for the
trigonometry is, more or less, artificial. Hence we leave the
following problem. \vskip 0.5pc \noindent{\bf Problem 1.} Find a
natural proof for the trigonometry on the extended hyperbolic space
or the extended de Sitter space. \vskip 0.5pc

 We can
consider a triangle area formula on the extended hyperbolic space
and the extended de Sitter space.  Particularly, an  area formula
for a triangle with angles $A,B,C$ is represented by $\pi-A-B-C$ on
$\Bbb S^n_H$ (naturally $A+B+C-\pi$ on $\Bbb S^n_S$) (see \cite{2}).

If we apply the cosine law on the extended hyperbolic space to the
area formula $S=\pi-A-B-C$, then we obtain another area formula
$S_1(a,b,c)$ with three edge length variables $a,b,c$,
$$\al S_1=\pi-\cos^{-1}\left(\frac{\cosh b\cosh c-\cosh a}{\sinh b \sinh c}\right)-\cos^{-1}\left(\frac{\cosh a\cosh c-\cosh b}{\sinh a \sinh
c}\right)\\
-\cos^{-1}\left(\frac{\cosh a\cosh b-\cosh c}{\sinh a \sinh
b}\right).\eal$$

We already know another area formula $S_2(a,b,c)$ (see \cite{1},
there is a misprint that is easily checked by considering
$a=b=c=\infty$) for a triangle on the hyperbolic space,
$$\tan^2 \frac{S_2}{4}=\tanh \frac{p}{2} \tanh \frac{p-a}{2} \tanh \frac{p-b}{2} \tanh \frac{p-c}{2},\quad \text{where }p=\frac{a+b+c}{2}. $$

Two function $S_1$ and $S_2$ are complex multi-valued functions on
$\Bbb C^3$, and coincide each other when the triangle lies on the
hyperbolic space, i.e., have the same value on a domain $U\subset
\Bbb R^3\subset\Bbb C^3$. Hence $S_1$ and $S_2$ coincide each other
on $\Bbb C^3$ with the same branch cuts, and so they have the same
value for a triangle on the extended hyperbolic space.

We know the principle $i\vv v\vv_S=\vv v \vv_H$, which comes from
Convention \ref{1.2} and \ref{1.4}, for two norms $\vv v\vv_H$ and
$\vv v \vv_S$ of any direction tangent vector $v$, where $\vv
\cdot\vv_H$ (resp. $\vv\cdot\vv_S$) denotes a vector norm on the
extended hyperbolic space (resp. extended de Sitter space). So any
2-dimensional volume elements $dV_H$ and  $dV_S$ for a given point
on $\Bbb S^n_H$ and $\Bbb S^n_S$, respectively, give a natural
relation $i^2\cdot dV_S=dV_H$.

Therefore we get an area formula for a triangle on the extended de
sitter space as well as on the spherical space,
$$\tan^2 \frac{-S_2}{4}=\tanh \frac{pi}{2} \tanh \frac{(p-a)i}{2} \tanh \frac{(p-b)i}{2} \tanh \frac{(p-c)i}{2},$$
$$\text{i.e., }\quad \tan^2 \frac{S_2}{4}=\tan \frac{p}{2} \tan \frac{p-a}{2} \tan \frac{p-b}{2} \tan \frac{p-c}{2},
\quad \text{where }p=\frac{a+b+c}{2}.$$

As a result, we can anticipate the following principle by Cho and
Kim. \vskip 0.5pc \noindent{\bf Problem 2.} If an analytic
(multi-valued) formula with geometric quantity variables is
satisfied on the hyperbolic space, then we can obtain the
corresponding formula on the spherical space by changing of all
variables with a {\bf principle that $k$-dimensional hyperbolic
variable is replaced by $i^k\times$ corresponding $k$-dimensional
spherical variable}, for example, hyperbolic angle $\theta$
$\rightarrow$ spherical angle $\theta$ and hyperbolic length $l$
$\rightarrow$ $i\times$ spherical length $l$ and so on. \vskip 0.5pc

In fact, if an $n$-dimensional (the highest dimension among the
variables' dimensions) analytic formula is satisfied on the
$n$-dimensional hyperbolic space (resp. spherical space) and if we
prove that the analytic formula also holds in the
$(n+1)$-dimensional extended hyperbolic space (resp. extended de
Sitter space), then Problem 2 is automatically satisfied by the
comparison of the extended hyperbolic space and the extended de
Sitter space.

In the above problem, if we change the contour for the extended
space, the value $i^k$ can be replaced by $(-i)^k$,
$k=1,2,3,\ldots$. So all analytic formula must have a symmetry for
$i$ and $-i$, i.e.,
$$f(\ldots,i^k\cdot k\text{-dim var.},\ldots,i^n\cdot n\text{-dim var.})=
f(\ldots,(-i)^k\cdot k\text{-dim var.},\ldots,(-i)^n\cdot
n\text{-dim var.}).$$

For hyperbolic and spherical triangles, Lobachevsky even knew the
principle for the hyperbolic and spherical trigonometry laws.

For $n$-dimensional hyperbolic and spherical simplices, Vinberg
\cite{V} clarified the principle for the $n$-dimensional volume of
the simplex and its dihedral angles.

\begin{flushleft}
{\sc  Department of Mathematics, University of Seoul, Seoul 130-743,
Korea\\[5pt]
  E-mail: yhcho@uos.ac.kr\\
}
\end{flushleft}
\end{document}